\documentclass[twoside,11pt]{article}

\usepackage{jmlr2e}

\usepackage{amsmath}
\usepackage{amsfonts,dsfont} 

\newcommand{\latin}[1]{\emph{#1}}

\newtheorem{proc}{Algorithm}

\newcommand{\egaldef}{:=} 
\newcommand{\defegal}{=:} 
\newcommand{\flens}{\mapsto} 
\newcommand{\flapp}{\mapsto} 
\newcommand{\telque}{\, \mbox{ s.t. } \,} 
\newcommand{\grandO}{\ensuremath{\mathcal{O}}}

\newcommand{\1}{\mathds{1}} 

\newcommand{\R}{\mathbb{R}} 
 
\newcommand{\X}{\mathcal{X}}
\newcommand{\Y}{\mathcal{Y}}

\DeclareMathOperator{\card}{Card} 

\newcommand{\mini}[2]{#1 \wedge #2}

\newcommand{\maxi}[2]{#1 \vee #2}

\newcommand{\paren}[1]{\left( \left. #1 \right. \right)} 
\newcommand{\croch}[1]{\left[ \left. #1 \right. \right]} 
\newcommand{\set}[1]{\left\{ \left. #1 \right. \right\}}
\newcommand{\absj}[1]{\left\lvert #1 \right\rvert} 
\providecommand{\norm}[1]{\left \lVert #1 \right\rVert}
\newcommand{\carre}[1]{\left(#1\right)^2}

\newcommand{\Prob}{\mathbb{P}} 
\newcommand{\E}{\mathbb{E}} 
\DeclareMathOperator{\var}{var} 
\newcommand{\sachant}{\, \right| \left. \,} 
 
\newcommand{\loi}{\mathcal{L}} 
\DeclareMathOperator{\Leb}{Leb} 

\newcommand{\bayes}{s}
\newcommand{\perte}[1]{\ell\paren{\bayes , #1 }}
\newcommand{\pertedeux}[2]{\ell\paren{ #1 , #2 }}
\newcommand{\ERM}{\widehat{s}}

\newcommand{\M}{\mathcal{M}}
\newcommand{\mM}{m \in \M}

\newcommand{\mo}{m^{\star}} 
\newcommand{\mh}{\widehat{m}}

\DeclareMathOperator{\pen}{pen}
\DeclareMathOperator{\crit}{crit}
\newcommand{\penid}{\pen_{\mathrm{id}}} 

\newcommand{\perteempc}[1]{\overline{\gamma}_n(#1)}

\newcommand{\esm}{\varepsilon_{\star,m}} 

\newcommand{\hyp}[1]{\ensuremath{\mathbf{(#1)}}} 
\newcommand{\hypAb}{\hyp{Ab}} 
\newcommand{\hypAn}{\hyp{An}} 
\newcommand{\hypAp}{\hyp{Ap}} 
\newcommand{\hypArXl}{\hyp{Ar^{X}_{\ell}}} 
\newcommand{\hypPpoly}{\hyp{P1}} 
\newcommand{\hypPrich}{\hyp{P2}} 
\newcommand{\hypPrichMin}{\hyp{P2+}} 

\newcommand{\hypHP}{\hyp{SH5}}
\newcommand{\hypHPm}{\hyp{SH2}}
\newcommand{\hypApu}{\hyp{Ap_u}} 

\newcommand{\crXl}{\ensuremath{c_{\mathrm{r},\ell}^X}}
\newcommand{\cbiasmaj}{C_{+}} 
\newcommand{\cbiasmin}{C_{-}} 
\newcommand{\betamaj}{\beta_{+}} 
\newcommand{\betamin}{\beta_{-}} 
\newcommand{\aM}{\alpha_{\M}}
\newcommand{\cM}{c_{\M}}
\newcommand{\sigmin}{\sigma_{\min}} 
\newcommand{\Qmp}{\ensuremath{Q_{m}^{(p)}}} 

\newcommand{\Kminpr}{K_1} 
\newcommand{\Kmindim}{K_2} 
\newcommand{\Koptpr}{K_3} 
\newcommand{\Kdimjump}{K_4} 
\newcommand{\Kdimjumpd}{K_5} 

\newcommand{\delc}{\overline{\delta}} 
\newcommand{\punmin}{\widetilde{p_1}} 

\newcommand{\Il}{I_{\lambda}}
\newcommand{\lamm}{\lambda \in \Lambda_m} 

\newcommand{\pl}{p_{\lambda}} 
\newcommand{\betl}{\beta_{\lambda}} 
\newcommand{\sigl}{\sigma_{\lambda}}

\newcommand{\tl}{t_{\lambda}}
\newcommand{\ul}{u_{\lambda}}

\newcommand{\phl}{\widehat{p}_{\lambda}} 
\newcommand{\bethl}{\widehat{\beta}_{\lambda}} 

\jmlrheading{0}{0000}{0}{3/08; Revised 9/08}{0/00}{Sylvain Arlot and Pascal Massart}

\ShortHeadings{Data-driven Calibration of Penalties}{Arlot and Massart}
\firstpageno{1}

\begin{document}

\title{Data-driven Calibration of Penalties for Least-Squares Regression}

\author{\name Sylvain Arlot \email sylvain.arlot@math.u-psud.fr \\
       \name Pascal Massart \email pascal.massart@math.u-psud.fr \\
       \addr Univ Paris-Sud, UMR 8628 \\Laboratoire de Mathematiques\\ Orsay, F-91405 ;
CNRS, Orsay, F-91405 ;\\
INRIA Saclay, Projet Select}

\editor{John Lafferty}

\maketitle

\begin{abstract}
Penalization procedures often suffer from their dependence on multiplying factors, whose optimal values are either unknown or hard to estimate from data. We propose a completely data-driven calibration algorithm for these parameters in the least-squares regression framework, without assuming a particular shape for the penalty. 
Our algorithm relies on the concept of minimal penalty, recently introduced by Birg\'{e} and Massart (2007) in the context of penalized least squares for Gaussian homoscedastic regression. 
On the positive side, the minimal penalty can be evaluated from the data themselves, leading to a data-driven estimation of an optimal penalty which can be used in practice; on the negative side, their approach heavily relies on the homoscedastic Gaussian nature of their stochastic framework. 

The purpose of this paper is twofold: stating a more general heuristics for designing a
data-driven penalty (the \textit{slope heuristics}) and proving that it works for penalized least-squares regression with a random design, even for heteroscedastic non-Gaussian data. 
For technical reasons, some exact mathematical results will be proved only for regressogram bin-width selection. This is at least a first step towards further results, since the approach and the method that we use are indeed general.
\end{abstract}


\section{Introduction} \label{sec.intro}
In the last decades, model selection has received much interest, commonly through penalization. 
In short, penalization chooses the model minimizing the sum of the empirical risk (how well the algorithm fits the data) and of some measure of complexity of the model (called penalty); see FPE \citep{Aka:1969}, AIC \citep{Aka:1973}, Mallows' $C_p$ or $C_L$ \citep{Mal:1973}. 
Many other penalization procedures have been proposed since, among which Rademacher complexities \citep{Kol:2001,Bar_Bou_Lug:2002}, local Rademacher complexities \citep{Bar_Bou_Men:2005,Kol:2006}, bootstrap penalties \citep{Efr:1983}, resampling and $V$-fold penalties \citep{Arl:2008a,Arl:2008b}.

Model selection can target two different goals. 
On the one hand, a procedure is {\em efficient} (or asymptotically optimal) when its quadratic risk is asymptotically equivalent to the risk of the oracle. 
On the other hand, a procedure is {\em consistent} when it chooses the smallest true model asymptotically with probability one.
This paper deals with {\em efficient} procedures, without assuming the existence of a true model.

\medskip

A huge amount of literature exists about efficiency. First Mallows' $C_p$, Akaike's FPE and AIC are asymptotically optimal, as proved by \cite{Shi:1981} for Gaussian errors, by \cite{KCLi:1987} under suitable moment assumptions on the errors, and by \cite{Pol_Tsy:1990} under sharper moment conditions, in the Fourier case.
Non-asymptotic oracle inequalities (with some leading constant $C>1$) have been obtained by \cite{Bar_Bir_Mas:1999} and by \cite{Bir_Mas:2002} in the Gaussian case, and by \cite{Bar:2000,Bar:2002} under some moment assumptions on the errors. 
In the Gaussian case, non-asymptotic oracle inequalities with leading constant $C_n$ tending to 1 when $n$ tends to infinity have been obtained by \cite{Bir_Mas:2006}.

However, from the practical point of view, both AIC and Mallows' $C_p$ still present serious drawbacks. 
On the one hand, AIC relies on a strong asymptotic assumption, so that for small sample sizes, the optimal multiplying factor can be quite different from one. 
Therefore, corrected versions of AIC have been proposed \citep{Sug:1978,Hur_Tsa:1989}.
On the other hand, the optimal calibration of Mallows' $C_p$ requires the knowledge of the noise level $\sigma^2$, assumed to be constant. 
When real data are involved, $\sigma^2$ has to be estimated separately and independently from any model, which is a difficult task. 
Moreover, the best estimator of $\sigma^2$ (say, with respect to the quadratic error) quite unlikely leads to the most efficient model selection procedure. 
Contrary to Mallows' $C_p$, the data-dependent calibration rule defined in this article is not a ``plug-in'' method; it focuses directly on efficiency, which can improve significantly the performance of the model selection procedure.

Existing penalization procedures present similar or stronger drawbacks than AIC and Mallows' $C_p$, often because of a gap between theory and practice. 
For instance, oracle inequalities have only been proved for (global) Rademacher penalties multiplied by a factor two \citep{Kol:2001}, while they are used without this factor \citep{Loz:2000}. 
As proved by \citet[Chapter~9]{Arl:2007:phd}, this factor is necessary in general. 
Therefore, the optimal calibration of these penalties is really an issue.
The calibration problem is even harder for local Rademacher complexities: theoretical results hold only with large calibration constants, particularly the multiplying factor, and no optimal values are known. 
One of the purposes of this paper is to address the issue of optimizing the multiplying factor for general-shape penalties.

\medskip

Few automatic calibration algorithms are available. The most popular ones are certainly cross-validation methods \citep{All:1974,Sto:1974}, in particular $V$-fold cross-validation \citep{Gei:1975}, because these are general-purpose methods, relying on a widely valid heuristics. 
However, their computational cost can be high.
For instance, $V$-fold cross-validation requires the entire model selection procedure to be performed $V$ times for each candidate value of the constant to be calibrated. 
For penalties proportional to the dimension of the models, such as Mallows' $C_p$, alternative calibration procedures have been proposed by \cite{Geo_Fos:2000} and by \cite{She_Ye:2002}. 

A completely different approach has been proposed by \cite{Bir_Mas:2006} for calibrating dimensionality-based penalties. 
Since this article extends their approach to a much wider range of applications, let us briefly recall their main results.
In Gaussian homoscedastic regression with a fixed design, assume that each model is a finite-dimensional vector space. 
Consider the penalty $\pen(m) = K D_m$, where $D_m$ is the dimension of the model $m$ and $K>0$ is a positive constant, to be calibrated. 
First, there exists a {\em minimal} constant $K_{\min}$, such that the ratio between the quadratic risk of the chosen estimator and the quadratic risk of the oracle is asymptotically infinite if $K<K_{\min}$, and finite if $K>K_{\min}$.
Second, when $K = K^{\star} \egaldef 2 K_{\min}$, the penalty $K D_m$ yields an efficient model selection procedure.
In other words, {\em the optimal penalty is twice the minimal penalty}.
This relationship characterizes the ``slope heuristics'' of \cite{Bir_Mas:2006}.

A crucial fact is that the minimal constant $K_{\min}$ can be estimated from the data, since large models are selected if and only if $K < K_{\min}$. 
This leads to the following strategy for choosing $K$ from the data.
For every $K \geq 0$, let $\mh(K)$ be the model selected by minimizing the empirical risk penalized by $\pen(D_m) = K D_m$. 
First, compute $K_{\min}$ such that $D_{\mh(K)}$ is ``huge'' for $K < K_{\min}$ and ``reasonably small'' when $K \geq K_{\min}$; explicit values for ``huge'' and ``small'' are proposed in Section~\ref{HP.sec.algo.defKmin}. 
Second, define $\mh \egaldef \mh (2 K_{\min} )$.
Such a method has been successfully applied for multiple change points detection by  \cite{Leb:2005}. 

From the theoretical point of view, the issue for understanding and validating this approach is the existence of a minimal penalty. 
This question has been addressed for Gaussian homoscedastic regression with a fixed design by \cite{Bir_Mas:2002,Bir_Mas:2006} when the variance is known, and by \cite{Bar_Gir_Hue:2007} when the variance is unknown. 
Non-Gaussian or heteroscedastic data have never been considered. 
This article contributes to fill this gap in the theoretical understanding of penalization procedures.

\medskip

The calibration algorithm proposed in this article relies on a generalization of Birg\'e and Massart's slope heuristics (Section~\ref{HP.sec.cadre.heuristic}).
In Section~\ref{HP.sec.pratique}, the algorithm is defined in the least-squares regression framework, for general-shape penalties.
The shape of the penalty itself can be estimated from the data, as explained in Section~\ref{HP.sec.pratique.shape}.

The theoretical validation of the algorithm is provided in Section~\ref{HP.sec.theorie}, from the {\em non-asymptotic point of view}. 
Non-asymptotic means in particular that the collection of models is allowed to depend on $n$: in practice, it is usual to allow the number of explanatory variables to increase with the number of observations.
Considering models with a large number of parameters (for example of the order of a power of the sample size $n$) is also necessary to approximate functions belonging to a general approximation space. 
Thus, the non-asymptotic point of view allows us not to assume that the regression function is described with a small number of parameters. 

The existence of minimal penalties for {\em heteroscedatic regression with a random design} (Theorem~\ref{HP.the.mini}) is proved in Section~\ref{HP.sec.theorie.mini}. 
In Section~\ref{HP.sec.theorie.opt}, by proving that twice the minimal penalty has some optimality properties (Theorem~\ref{HP.the.opt}), we extend the so-called slope heuristics to heteroscedatic regression with a random design. 
Moreover, neither Theorem~\ref{HP.the.mini} nor Theorem~\ref{HP.the.opt} assume the data to be Gaussian; only mild moment assumptions are required.

For proving Theorems~\ref{HP.the.mini} and~\ref{HP.the.opt}, each model is assumed to be the vector space of piecewise constant functions on some partition of the feature space. 
This is indeed a restriction, but we conjecture that it is mainly technical, and that the slope heuristics remains valid at least in the general least-squares regression framework. 
We provide some evidence for this by proving two key concentration inequalities without the restriction to piecewise constant functions.
Another argument supporting this conjecture is that recently several simulation studies have shown that the slope heuristics can be used in several frameworks: mixture models \citep{Mau_Mic:2008}, clustering \citep{Bau:2007}, spatial statistics \citep{Ver:2008:phd}, estimation of oil reserves \citep{Lep:2002} and genomics \citep{Vil:2007:these}.
Although the slope heuristics has not been formally validated in these frameworks, this article is a first step towards such a validation, by proving that the slope heuristics can be applied whatever the shape of the ideal penalty. 

\medskip

This paper is organized as follows. The framework and the slope heuristics are described in Section~\ref{HP.sec.cadre}. The resulting algorithm is defined in Section~\ref{HP.sec.pratique}. The main theoretical results are stated in Section~\ref{HP.sec.theorie}. All the proofs are given in Appendix~\ref{HP.sec.proofs}.

\section{Framework} \label{HP.sec.cadre}
In this section, we describe the framework and the general slope heuristics.
\subsection{Least-squares regression} \label{HP.sec.regression}
Suppose we observe some data $(X_1,Y_1), \ldots (X_n,Y_n) \in \X \times \R$, independent with common distribution $P$, where the feature space $\X$ is typically a compact set of $\R^d$. 
The goal is to predict $Y$ given $X$, where $(X,Y) \sim P$ is a new data point independent of $(X_i,Y_i)_{1 \leq i \leq n}$. 
Denoting by $\bayes$ the regression function, that is $\bayes(x) = \E\croch{Y \sachant X = x}$ for every $x \in \X$, we can write 
\begin{equation} \label{HP.eq.donnees.reg} Y_i = \bayes(X_i) + \sigma(X_i) \epsilon_i \end{equation} 
where $\sigma: \X \mapsto \R$ is the heteroscedastic noise level and $\epsilon_i$ are i.i.d. centered noise terms, possibly dependent on $X_i$, but with mean 0 and variance 1 conditionally to $X_i$. 

The quality of a predictor $t: \X \mapsto \Y$ is measured by the (quadratic) prediction loss
\[ \E_{(X,Y) \sim P} \croch{ \gamma(t,(X,Y)) } \defegal P \gamma(t) \qquad \mbox{where} \quad \gamma(t,(x,y)) = \paren{ t(x) - y }^2 \]
is the least-squares contrast. 
The minimizer of $P \gamma(t)$ over the set of all predictors, called Bayes predictor, is the regression function $\bayes$. 
Therefore, the excess loss is defined as  
\[ \perte{t} \egaldef P\gamma\paren{t} - P\gamma\paren{\bayes} =  \E_{(X,Y)\sim P} \carre{t(X)-\bayes(X)} \enspace . \]
Given a particular set of predictors $S_m$ (called a {\em model}), we define the best predictor over $S_m$ as \[ \bayes_m \egaldef \arg\min_{t \in S_m} \set{ P \gamma(t) } \enspace , \]
with its empirical counterpart 
\[ \ERM_m \egaldef \arg\min_{t \in S_m} \set{ P_n \gamma(t) } \] (when it exists and is unique), where $P_n = n^{-1} \sum_{i=1}^n \delta_{(X_i,Y_i)}$. 
This estimator is the well-known 
{\em empirical risk minimizer}, also called least-squares estimator since $\gamma$ is the least-squares contrast.

\subsection{Ideal model selection} \label{HP.sec.cadre.mod_selec}
Let us assume that we are given a family of models $(S_m)_{\mM_n}$, hence a family of estimators $(\ERM_m)_{\mM_n}$ obtained by empirical risk minimization.
The model selection problem consists in looking for some data-dependent $\mh \in \M_n$ such that $\perte{\ERM_{\mh}}$ is as small as possible.
For instance, it would be convenient to prove some oracle inequality of the form
\[ \perte{\ERM_{\mh}} \leq C \inf_{\mM_n} \set{\perte{\ERM_m}} + R_n \]
in expectation or on an event of large probability, with leading constant $C$ close to 1 and $R_n = o(n^{-1})$.

General penalization procedures can be described as follows.
Let $\pen: \M_n \flapp \R^+$ be some penalty function, possibly data-dependent, and define  \begin{equation} \label{HP.eq.penalization} \mh \in \arg\min_{\mM_n} \set{\crit(m)} \qquad \mbox{with} \qquad \crit(m) \egaldef P_n \gamma(\ERM_m) + \pen(m) \enspace . \end{equation} 
Since the ideal criterion $\crit(m)$ is the true prediction error $P \gamma\paren{\ERM_m}$, the {\em ideal penalty} is 
\[ \penid(m) \egaldef P \gamma(\ERM_m) - P_n \gamma(\ERM_m) \enspace . \]
This quantity is unknown because it depends on the true distribution $P$. 
A natural idea is to choose $\pen(m)$ as close as possible to $\penid(m)$ for every $\mM_n$. 
We will show below, in a general setting, that when $\pen$ is a good estimator of the ideal penalty $\penid$, then $\mh$ satisfies an oracle inequality with leading constant $C$ close to 1.

\medskip

By definition of $\mh$,
\[ \forall \mM_n, \quad P_n \gamma(\ERM_{\mh}) \leq 
P_n \gamma(\ERM_m) + \pen(m) - \pen(\mh) \enspace .\]
For every $\mM_n$, we define
\begin{align*}
p_1(m) = P \paren{ \gamma(\ERM_m) - \gamma(\bayes_m) } \quad
 p_2(m) = P_n \paren{ \gamma(\bayes_m) - \gamma(\ERM_m) }  \quad
 \delta(m) =  (P_n - P) \paren{ \gamma(\bayes_m)}
\end{align*}
so that 
\begin{gather*}
\penid(m) = p_1(m) + p_2(m) - \delta(m) \\
\mbox{and} \qquad \perte{\ERM_{m}} = P_n \gamma(\ERM_{m}) + p_1(m) + p_2(m) - \delta(m) -
P\gamma(s) \enspace .\end{gather*}
Hence, for every $\mM_n$,
\begin{gather}\label{HP.eq.etape1} \perte{\ERM_{\mh}} + (\pen - \penid) (\mh ) \leq \perte{\ERM_m} +  (\pen - \penid) (m) \enspace . 
\end{gather}
Therefore, in order to derive an oracle inequality from \eqref{HP.eq.etape1}, it is sufficient to show that for every $\mM_n$, $\pen(m)$ is close to $\penid(m)$.

\subsection{The slope heuristics} \label{HP.sec.cadre.heuristic}
If the penalty is too big, the left-hand side of \eqref{HP.eq.etape1} is larger than $\perte{\ERM_{\mh}}$ so that \eqref{HP.eq.etape1} implies an oracle inequality, possibly with large leading constant $C$.
On the contrary, if the penalty is too small, the left-hand side of \eqref{HP.eq.etape1} may become negligible with respect to $\perte{\ERM_{\mh}}$ (which would make $C$ explode) or---worse---may be nonpositive.
In the latter case, no oracle inequality may be derived from \eqref{HP.eq.etape1}. 
We shall see in the following that $\perte{\ERM_{\mh}}$ blows up if and only if the penalty is smaller than some ``minimal penalty''.

Let us consider first the case $\pen(m) = p_2(m)$ in \eqref{HP.eq.penalization}. 
Then, $\E\croch{\crit(m)} = \E\croch{P_n \gamma\paren{\bayes_m}} = P \gamma\paren{\bayes_m}$, so that $\mh$ approximately minimizes its bias. 
Therefore, $\mh$ is one of the more complex models, and the risk of $\ERM_{\mh}$ is large. 
Let us assume now that $\pen(m) = K p_2(m)$. 
If $0<K<1$, $\crit(m)$ is a decreasing function of the complexity of $m$, so that $\mh$ is again one of the more complex models. 
On the contrary, if $K>1$, $\crit(m)$ increases with the complexity of $m$ (at least for the largest models), so that $\mh$ has a small or medium complexity. 
This argument supports the conjecture that the ``minimal amount of penalty'' required for the model selection procedure to work is $p_2(m)$. 

\medskip

In many frameworks such as the one of Section~\ref{HP.sec.cadre.histos}, it turns out that 
\[ \forall \mM_n, \qquad p_1(m) \approx p_2(m) \enspace . \] 
Hence, the ideal penalty $\penid(m) \approx p_1(m) + p_2(m) $ is close to $2 p_2(m)$. 
Since $p_2(m)$ is a ``minimal penalty'', the optimal penalty is close to twice the minimal penalty:
\[ \penid(m) \approx 2 \pen_{\min}(m) \enspace . \]
This is the so-called ``slope heuristics'', first introduced by \cite{Bir_Mas:2006} in a Gaussian homoscedastic setting.
Note that a formal proof of the validity of the slope heuristics has only been given for Gaussian homoscedastic least-squares regression with a fixed design \citep{Bir_Mas:2006}; up to the best of our knowledge, the present paper yields the second theoretical result on the slope heuristics.

\medskip

This heuristics has some applications because the minimal penalty can be estimated from the data. 
Indeed, when the penalty smaller than $\pen_{\min}$, the selected model $\mh$ is among the more complex. 
On the contrary, when the penalty is larger than $\pen_{\min}$, the complexity of $\mh$ is much smaller. 
This leads to the algorithm described in the next section.

\section{A data-driven calibration algorithm}
\label{HP.sec.pratique}
Now, a data-driven calibration algorithm for penalization procedures can be defined, generalizing a method proposed by \cite{Bir_Mas:2006} and implemented by \cite{Leb:2005}.

\subsection{The general algorithm} \label{HP.sec.pratique.data-driven}
Assume that the shape $\pen_{\mathrm{shape}}: \M_n \flens \R^+$ of the ideal penalty is known, from some prior knowledge or because it had first been estimated, see Section~\ref{HP.sec.pratique.shape}. 
Then, the penalty $K^{\star} \pen_{\mathrm{shape}}$ provides an approximately optimal procedure, for some unknown constant $K^{\star}>0$. 
The goal is to find some $\widehat{K}$ such that $\widehat{K} \pen_{\mathrm{shape}}$ is approximately optimal. 

Let $D_m$ be some known complexity measure of the model $\mM_n$.
Typically, when the models are finite-dimensional vector spaces, $D_m$ is the dimension of $S_m$.
According to the ``slope heuristics'' detailed in Section~\ref{HP.sec.cadre.heuristic}, the following algorithm provides an optimal calibration of the penalty $\pen_{\mathrm{shape}}$.

\begin{proc}[Data-driven penalization with slope heuristics] \label{HP.def.proc.gal.pente} 
\hfill
\begin{enumerate}
\item Compute the selected model $\mh(K)$ as a function of $K>0$
\begin{equation*} 
\mh(K) \in \arg\min_{\mM_n} \left\{ P_n \gamma(\ERM_m) + K \pen_{\mathrm{shape}} (m) \right\} \enspace .
\end{equation*}
\item Find $\widehat{K}_{\min} > 0$ such that $D_{\mh(K)}$ is ``huge'' for $K < \widehat{K}_{\min}$ and ``reasonably small'' for $K > \widehat{K}_{\min}$.
\item Select the model $\mh \egaldef \mh\paren{2\widehat{K}_{\min}}$.
\end{enumerate}
\end{proc}

A computationally efficient way to perform the first step of Algorithm~\ref{HP.def.proc.gal.pente} is provided in Section~\ref{HP.sec.algo.computKmin}.
The accurate definition of $\widehat{K}_{\min}$ is discussed in Section~\ref{HP.sec.algo.defKmin}, including explicit values for ``huge'' and ``reasonably small'').
Then, once $P_n \gamma\paren{\ERM_m}$ and $\pen_{\mathrm{shape}} (m)$ are known for every $\mM_n$, the complexity of Algorithm~\ref{HP.def.proc.gal.pente} is $\mathcal{O}(\card(\M_n)^2)$ (see Algorithm~\ref{HP.algo.step} and Proposition~\ref{HP.pro.algo.step}). 
This can be a decisive advantage compared to cross-validation methods, as discussed in Section~\ref{HP.sec.theorie.comment.cv}.

\subsection{Computation of $\paren{\mh(K)}_{K\geq 0}$} \label{HP.sec.algo.computKmin}
Step~1 of Algorithm~\ref{HP.def.proc.gal.pente} requires to compute $\mh(K)$ for every $K \in (0,+\infty)$. 
A computationally efficient way to perform this step is described in this subsection. 

We start with some notations: 
\[ \forall \mM_n, \qquad f(m) = P_n \gamma\paren{\ERM_m} \qquad
\qquad g(m) = \pen_{\mathrm{shape}} (m) \]
\[ \mbox{and} \qquad \forall K\geq 0, \qquad \mh(K) \in 
\arg\min_{\mM_n} \set{ f(m) + K g(m) } \enspace . \]
Since the latter definition can be ambiguous, let us choose any total ordering $\preceq$ on $\M_n$ such that $g$ is non-decreasing, which is always possible if $\M_n$ is at most countable. Then, $\mh(K)$ is defined as the smallest element of 
\[ E(K) \egaldef \arg\min_{\mM_n} \set{ f(m) + K g(m) }  \]
for $\preceq$.
The main reason why the whole trajectory $\paren{\mh(K)}_{K\geq 0}$ can be computed efficiently is its particular shape.

Indeed, the proof of Proposition~\ref{HP.pro.algo.step} shows that $K \flapp \mh(K)$ is piecewise constant, and non-increasing for $\preceq$. 
Then, the whole trajectory $\paren{\mh(K)}_{K\geq 0}$ can be summarized by
\begin{itemize}
\item the number of jumps $i_{\max} \in \set{0, \ldots, \card(\M_n) - 1}$,
\item the location of the jumps: an increasing sequence of nonnegative reals $(K_i)_{0 \leq i \leq i_{\max}+1}$, with $K_0=0$ and $K_{i_{\max}+1}=+\infty$, 
\item a non-increasing sequence of models $(m_i)_{0 \leq i \leq i_{\max}}$, 
\end{itemize}
\[ \mbox{with} \qquad \forall i \in \set{0, \ldots, i_{\max}}, \quad \forall K \in \left[ K_i , K_{i+1} \right) , \qquad \mh(K) = m_i \enspace . \]

\begin{proc}[Step~1 of Algorithm~\ref{HP.def.proc.gal.pente}] \label{HP.algo.step}
For every $\mM_n$, define $f(m) = P_n \gamma\paren{\ERM_m}$ and $g(m) = \pen_{\mathrm{shape}} (m)$. Choose $\preceq$ any total ordering on $\M_n$ such that $g$ is non-decreasing.
\begin{itemize}
\item Init: $K_0 \egaldef 0$, $m_0 \egaldef  \arg\min_{\mM_n} \set{ f(m) }$ (when this minimum is attained several times, $m_0$ is defined as the smallest one with respect to $\preceq$).
\item Step $i$, $i \geq 1$: 
Let 
\[ G(m_{i-1}) \egaldef \set{\mM_n \telque f(m) > f(m_{i-1}) \quad \mbox{ and } \quad g(m) < g(m_{i-1})} \enspace . \] If $G(m_{i-1})=\emptyset$, then put $K_i = + \infty$, $i_{\max}=i-1$ and stop. Otherwise,
\begin{gather} \label{HP.def.algo.step.Ki} 
K_{i} \egaldef \inf\set{\frac{f(m) - f(m_{i-1})} {g(m_{i-1}) - g(m)} \telque m \in G(m_{i-1}) } \\
\mbox{and} \qquad m_i \egaldef \min_{\preceq} F_i \qquad \mbox{with} \quad F_i \egaldef \arg\min_{m \in G(m_{i-1})} \set{\frac{f(m) - f(m_{i-1})} {g(m_{i-1}) - g(m)} } \enspace . \notag
\end{gather}
\end{itemize}
\end{proc}
\begin{proposition}[Correctness of Algorithm~\ref{HP.algo.step}] \label{HP.pro.algo.step} If $\M_n$ is finite, Algorithm~\ref{HP.algo.step} terminates and $i_{\max} \leq \card(\M_n) - 1$. 
With the notations of Algorithm~\ref{HP.algo.step}, let $\mh(K)$ be the smallest element of 
\[ E(K) \egaldef \arg\min_{\mM_n} \set{ f(m) + K g(m) } \qquad \mbox{with respect to } \preceq . \]
Then, $(K_i)_{0 \leq i \leq i_{\max}+1}$ is increasing and $\forall i \in \set{0, \ldots, i_{\max}-1}$, $\forall K \in [K_i, K_{i+1})$, $\mh(K) = m_i$.
\end{proposition}
It is proved in Section~\ref{HP.sec.algo.proof}. In the change-point detection framework, a similar result has been proved by \cite{Lav:2005}.

Proposition~\ref{HP.pro.algo.step} also gives an upper bound on the computational complexity of Algorithm~\ref{HP.algo.step};
since the complexity of each step is $\grandO(\card\M_n)$, Algorithm~\ref{HP.algo.step} requires less than $\grandO(i_{\max} \card \M_n) \leq \grandO((\card \M_n)^2)$ operations.
In general, this upper bound is pessimistic since $i_{\max} \ll \card \M_n$. 

\subsection{Definition of $\widehat{K}_{\min}$} \label{HP.sec.algo.defKmin}
Step~2 of Algorithm~\ref{HP.def.proc.gal.pente} estimates $\widehat{K}_{\min}$ such that $\widehat{K}_{\min} \pen_{\mathrm{shape}}$ is the minimal penalty. 
The purpose of this subsection is to define properly $\widehat{K}_{\min}$ as a function of $(\mh(K))_{K>0}$.

According to the slope heuristics described in Section~\ref{HP.sec.cadre.heuristic}, $\widehat{K}_{\min}$ corresponds to a 
``complexity jump''. 
If $K<\widehat{K}_{\min}$, $\mh(K)$ has a large complexity, whereas if $K>\widehat{K}_{\min}$, $\mh(K)$ has a small or medium complexity. 
Therefore, the two following definitions of $\widehat{K}_{\min}$ are natural.

Let $D_{\mathrm{thresh}}$ be the largest ``reasonably small'' complexity, meaning the models with larger complexities should not be selected. 
When $D_m$ is the dimension of $S_m$ as a vector space, $D_{\mathrm{thresh}} \propto n/(\ln(n))$ or $n/(\ln (n))^2$ are natural choices since the dimension of the oracle is likely to be of order $n^{\alpha}$ for some $\alpha \in (0,1)$. 
Then, define 
\begin{equation} \label{def.Kmin.thresh}
\tag{thresh} \widehat{K}_{\min} \egaldef  \inf\set{ K >0  \telque D_{\mh(K)} \leq D_{\mathrm{thresh}} } \enspace . 
\end{equation}
With this definition, Algorithm~\ref{HP.algo.step} can be stopped as soon as the threshold is reached. 

Another idea is that $\widehat{K}_{\min}$ should match with the largest complexity jump: 
\begin{equation}  \label{def.Kmin.maxjump}
\tag{max jump}  \widehat{K}_{\min} \egaldef K_{i_{\mathrm{jump}}} \quad \mbox{with} \quad i_{\mathrm{jump}} = \arg\max_{i \in \set{0, \ldots, i_{\max} - 1}} \set{ D_{m_{i+1}} - D_{m_{i}}} \enspace . 
\end{equation}
In order to ensure that there is a clear jump in the sequence $(D_{m_i})_{i \geq 0}$, it may be useful to add a few models of large complexity.

\medskip

\begin{figure}
\begin{minipage}[b]{.46\linewidth}
\includegraphics[width=0.95\textwidth]{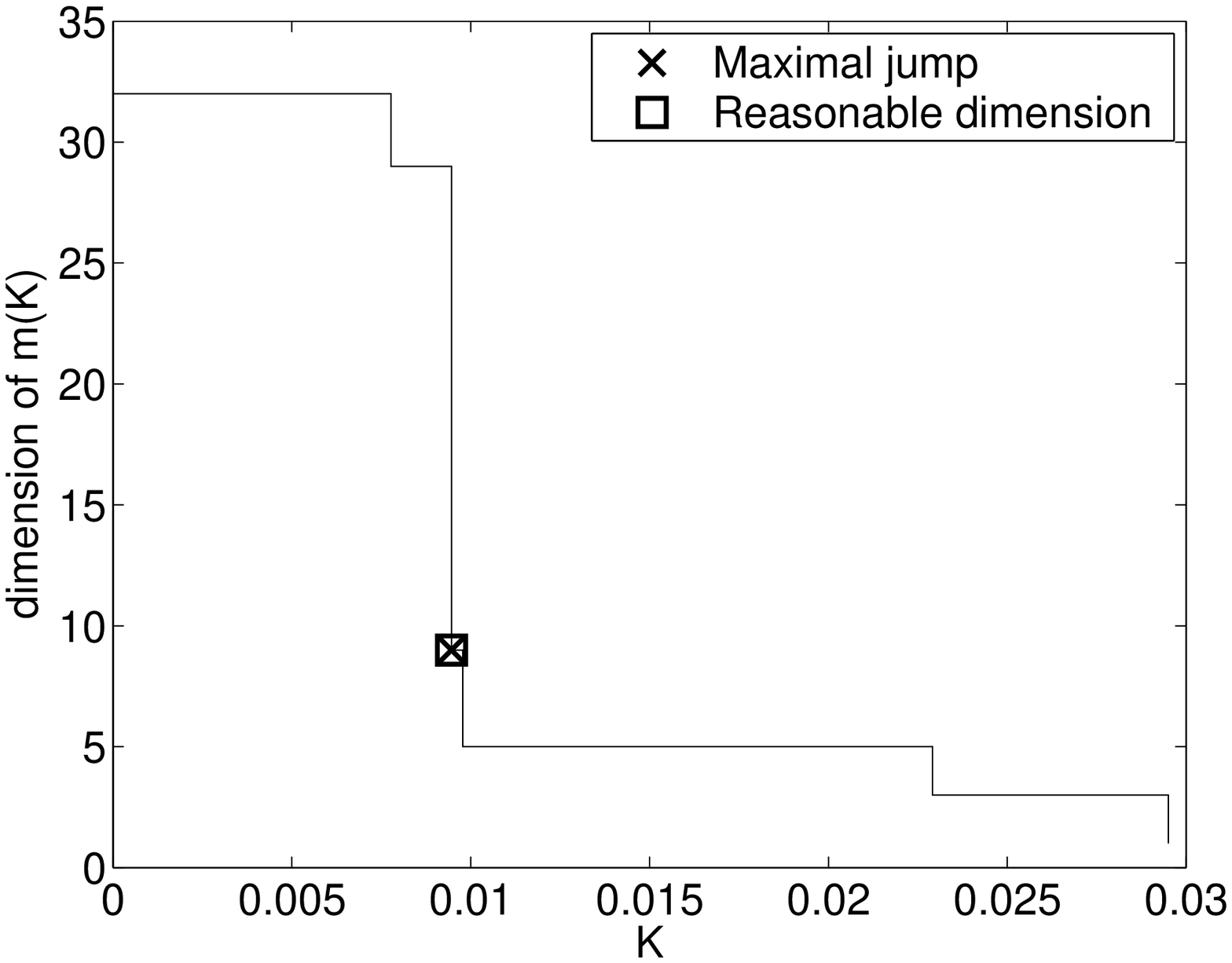}
\begin{center}
(a) One clear jump.
\end{center}
\end{minipage}
\begin{minipage}[b]{.46\linewidth}
\includegraphics[width=0.95\textwidth]{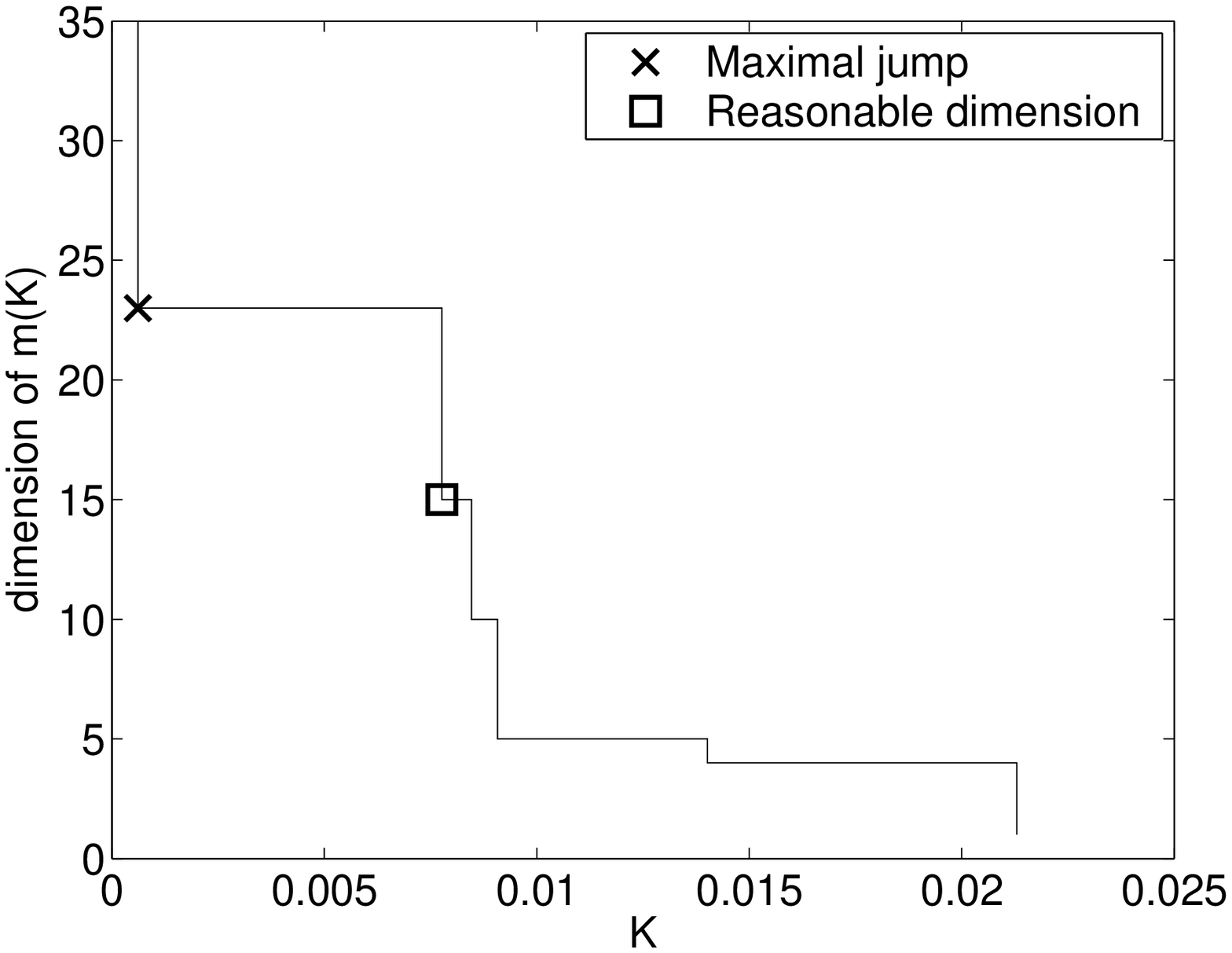}
\begin{center}
(b) Two jumps, two values for $\widehat{K}_{\min}$.
\end{center}
\end{minipage}
\caption{$D_{\mh(K)}$ as a function of $K$ for two different samples. 
Data are simulated according to \eqref{HP.eq.donnees.reg} with $n=200$, $X_i \sim\mathcal{U}([0,1])$, $\epsilon_i\sim\mathcal{N}(0,1)$, $\bayes(x)=\sin(\pi x)$ and $\sigma\equiv 1$. 
The models $(S_m)_{\mM_n}$ are the sets of piecewise constant functions on regular partitions of $[0,1]$, with dimensions between 1 and $n/(\ln(n))$.
The penalty shape is $\pen_{\mathrm{shape}}(m)= D_m$ and the dimension threshold is $D_{\mathrm{thresh}} = 19 \approx n/(2 \ln(n))$.
See experiment S1 by \citet[Section~6.1]{Arl:2008b} for details. \label{HP.fig.sauts}}
\end{figure}

As an illustration, we compared the two definitions above (``threshold'' and ``maximal jump'') on $1\,000$ simulated samples. 
The exact simulation framework is described below Figure~\ref{HP.fig.sauts}. 
Three cases occured:
\begin{enumerate}
\item There is one clear jump.
Both definitions give the same value for $\widehat{K}_{\min}$.
This occured for about $85\%$ of the samples; an example is given on Figure~\ref{HP.fig.sauts}a.
\item There are several jumps corresponding to close values of $K$. 
Definitions \eqref{def.Kmin.thresh} and \eqref{def.Kmin.maxjump} give slightly different values for $\widehat{K}_{\min}$, but the selected models $\mh\paren{2 \widehat{K}_{\min}}$ are equal.
This occured for about $8.5\%$ of the samples.
\item There are several jumps corresponding to distant values of $K$.
Definitions \eqref{def.Kmin.thresh} and \eqref{def.Kmin.maxjump} strongly disagree, giving different selected models $\mh\paren{2 \widehat{K}_{\min}}$ at final.
This occured for about $6.5\%$ of the samples; an example is given on Figure~\ref{HP.fig.sauts}b.
\end{enumerate}
The only problematic case is the third one, in which an arbitrary choice has to be made between definitions \eqref{def.Kmin.thresh} and \eqref{def.Kmin.maxjump}.

With the same simulated data, we have compared the prediction errors of the two methods by estimating the constant $C_{\mathrm{or}}$ that would appear in some oracle inequality,
\[ C_{\mathrm{or}} \egaldef \frac{\E\croch{ \perte{\ERM_{\mh}} }} {\E\croch{ \inf_{\mM_n} \set{ \perte{\ERM_m} } }} \enspace .\]
With definition \eqref{def.Kmin.thresh} $C_{\mathrm{or}} \approx 1.88$; with definition \eqref{def.Kmin.maxjump} $C_{\mathrm{or}} \approx 2.01$. 
For both methods, the standard error of the estimation is $0.04$.
As a comparison, Mallows' $C_p$ with a classical estimator of the variance $\sigma^2$ has an estimated performance $C_{\mathrm{or}} \approx 1.93$ on the same data. 

\medskip

The overall conclusion of this simulation experiment is that Algorithm~\ref{HP.def.proc.gal.pente} can be competitive with Mallows' $C_p$ in a framework where Mallows' $C_p$ is known to be optimal.
Definition \eqref{def.Kmin.thresh} for $\widehat{K}_{\min}$ seems slightly more efficient than \eqref{def.Kmin.maxjump}, but without convincing evidence.
Indeed, both definitions depend on some arbitrary choices: the value of the threshold $D_{\mathrm{thresh}}$ in \eqref{def.Kmin.thresh}, the maximal complexity among the collection of models $(S_m)_{\mM_n}$ in \eqref{def.Kmin.maxjump}. 
When $n$ is small, say $n=200$, choosing $D_{\mathrm{thresh}}$ is tricky since $n/(2\ln(n))$ and $\sqrt{n}$ are quite close. 
Then, the difference between \eqref{def.Kmin.thresh} and \eqref{def.Kmin.maxjump} is likely to come mainly from the particular choice $D_{\mathrm{thresh}} = 19$ than from basic differences between the two definitions.

In order to estimate $\widehat{K}_{\min}$ as automatically as possible, we suggest to combine the two definitions; when the selected models $\mh(2 \widehat{K}_{\min})$ differ, send a warning to the final user advising him to look at the curve $K \flapp D_{\mh(K)}$ himself; otherwise, remain confident in the automatic choice of $\mh(2 \widehat{K}_{\min})$.

\subsection{Penalty shape} \label{HP.sec.pratique.shape}
For using Algorithm~\ref{HP.def.proc.gal.pente} in practice, it is necessary to know \latin{a priori}, or at least to estimate, the optimal shape $\pen_{\mathrm{shape}}$ of the penalty. 
Let us explain how this can be achieved in different frameworks.

\medskip

The first example that comes to mind is $\pen_{\mathrm{shape}}(m) = D_m$. 
It is valid for homoscedastic least-squares regression on linear models, as shown by several papers mentioned in Section~\ref{sec.intro}. 
Indeed, when $\card(\M_n)$ is smaller than some power of $n$, Mallows' $C_p$ penalty---defined by $\pen(m) = 2 \E\croch{\sigma^2(X)} n^{-1} D_m$\,---is well known to be asymptotically optimal. 
For larger collections $\M_n$, more elaborate results \citep{Bir_Mas:2002,Bir_Mas:2006} have shown that a penalty proportional to $\ln(n) \E\croch{\sigma^2(X)} n^{-1} D_m$ and depending on the size of $\M_n$ is asymptotically optimal.

Algorithm~\ref{HP.def.proc.gal.pente} then provides an alternative to plugging an estimator of $\E\croch{\sigma^2(X)}$ into the above penalties. 
Let us detail two main advantages of our approach.
First, we avoid the difficult task of estimating $\E\croch{\sigma^2(X)}$ without knowing in advance some model to which the true regression function belongs.
Algorithm~\ref{HP.def.proc.gal.pente} provides a model-free estimation of the factor multiplying the penalty.
Second, the estimator $\widehat{\sigma^2}$ of $\E\croch{\sigma^2(X)}$ with the smallest quadratic risk is certainly far from being the optimal one for model selection. 
For instance, underestimating the multiplicative factor 
is well-known to lead to poor performances, whereas overestimating the multiplicative factor 
does not increase much the prediction error in general.
Then, a good estimator of $\E\croch{\sigma^2(X)}$ for model selection should overestimate it with a probability larger than $1/2$. 
Algorithm~\ref{HP.def.proc.gal.pente} satisfies this property automatically because $\widehat{K}_{\min}$ so that  the selected model cannot be too large.

In short, {\em Algorithm~\ref{HP.def.proc.gal.pente} with $\pen_{\mathrm{shape}}(m) = D_m$ is quite different from a simple plug-in version of Mallows' $C_p$}. It leads to a really {\em data-dependent penalty, which may perform better in practice than the best deterministic penalty $K^{\star} D_m$}.

\medskip

In a more general framework, Algorithm~\ref{HP.def.proc.gal.pente} allows to choose a different shape of penalty $\pen_{\mathrm{shape}}$. 
For instance, in the heteroscedastic least-squares regression framework of Section~\ref{HP.sec.regression}, the optimal penalty is no longer proportional to the dimension $D_m$ of the models. This can be shown from computations made by \cite[Proposition~1]{Arl:2008b} when $S_m$ is assumed to be the vector space of piecewise constant functions on a partition $\paren{\Il}_{\lamm}$ of $\X$:
\begin{equation} \label{HP.eq.Epenid}
\E\croch{\penid(m)} = \E \croch{ (P-P_n)\gamma\paren{\ERM_m} } \approx \frac{2}{n} \sum_{\lamm} \E\croch{ \sigma(X)^2 \sachant X \in \Il } \enspace .
\end{equation}
An exact result has been proved by \citet[Proposition~1]{Arl:2008b}. 
Moreover, \cite{Arl:2008:shape} gave an example of model selection problem in which no penalty proportional to $D_m$ can be asymptotically optimal.

A first way to estimate the shape of the penalty is simply to use \eqref{HP.eq.Epenid} to compute $\pen_{\mathrm{shape}}$, when both the distribution of $X$ and the shape of the noise level $\sigma$ are known. 
In practice, one has seldom such a prior knowledge.

We suggest in this situation to use {\em resampling penalties} \citep{Efr:1983,Arl:2008b}, or {\em $V$-fold penalties} \citep{Arl:2008a} which have much smaller computational costs. 
Up to a multiplicative factor (automatically estimated by Algorithm~\ref{HP.def.proc.gal.pente}), these penalties should estimate correctly $\E\croch{\penid(m)}$ in any framework. 
In particular, resampling and $V$-fold penalties are asymptotically optimal in the heteroscedastic least-squares regression framework \citep{Arl:2008a,Arl:2008b}.

\subsection{The general prediction framework} \label{HP.sec.pratique.gal}
Section~\ref{HP.sec.cadre} and definition of Algorithm~\ref{HP.def.proc.gal.pente} have restricted ourselves to the least-squares regression framework. 
Actually, this is not necessary at all to make Algorithm~\ref{HP.def.proc.gal.pente} well-defined, so that it can naturally be extended to the general prediction framework. 
More precisely, the $(X_i,Y_i)$ can be assumed to belong to $\X\times\Y$ for some general $\Y$, and $\gamma : S \times (\X \times \Y) \flens [0;+\infty)$ any contrast function. 
In particular, $\Y = \set{0,1}$ leads to the binary classification problem, for which a natural contrast function is the $0$--$1$ loss $\gamma(t;(x,y))=\1_{t(x) \neq y}$.
In this case, the shape of the penalty $\pen_{\mathrm{shape}}$ can for instance be estimated with the global or local Rademacher complexities mentioned in Section~\ref{sec.intro}.

However, a natural question is whether the slope heuristics of Section~\ref{HP.sec.cadre.heuristic}, upon which Algorithm~\ref{HP.def.proc.gal.pente} relies, can be extended to the general framework. 
Several concentration results used to prove the validity of the slope heuristics in the least-squares regression framework in this article are valid in a general setting including binary classification.
Even if the factor 2 coming from the closeness of $\E\croch{p_1}$ and $\E\croch{p_2}$ (see Section~\ref{HP.sec.cadre.heuristic}) may not be universally valid, we conjecture that Algorithm~\ref{HP.def.proc.gal.pente} can be used in other settings than least-squares regression.
Moreover, as already mentioned at the end of Section~\ref{sec.intro}, empirical studies have shown that Algorithm~\ref{HP.def.proc.gal.pente} can be successfully applied to several problems, with different shapes for the penalty.
To our knowledge, to give a formal proof of this fact remains an interesting open problem.

\section{Theoretical results} \label{HP.sec.theorie}
Algorithm~\ref{HP.def.proc.gal.pente} mainly relies on the ``slope heuristics'', developed in Section~\ref{HP.sec.cadre.mod_selec}. 
The goal of this section is to provide a theoretical justification of this heuristics.

It is split into two main results. 
First, Theorem~\ref{HP.the.mini} provides lower bounds on $D_{\mh}$ and the risk of $\ERM_{\mh}$ when the penalty is smaller than $\pen_{\min} (m) \egaldef \E\croch{p_2(m)}$. 
Second, Theorem~\ref{HP.the.opt} is an oracle inequality with leading constant almost one when $\pen (m) \approx 2\E\croch{p_2 (m)}$, relying on \eqref{HP.eq.etape1} and the comparison $p_1 \approx p_2$. 

\medskip

In order to prove both theorems, two probabilistic results are necessary. 
First, $p_1$, $p_2$ and $\delta$ concentrate around their expectations; for $p_2$ and $\delta$, it is proved in a general framework in Appendix~\ref{HP.sec.conc}. 
Second, $\E\croch{p_1(m)} \approx \E\croch{p_2(m)}$ for every $\mM_n$. 
The latter point is quite hard to prove in general, so that we must make an assumption on the models. 
Therefore, in this section, we restrict ourselves to the regressogram case, assuming that for every $\mM_n$, $S_m$ is the set of piecewise constant functions on some fixed partition $\paren{\Il}_{\lamm}$ of $\X$.
This framework is described precisely in the next subsection.
Although we do not consider regressograms as a final goal, the theoretical results proved for regressograms help to understand better how to use Algorithm~\ref{HP.def.proc.gal.pente} in practice.

\subsection{Regressograms} \label{HP.sec.cadre.histos}
Let $S_m$ be the the set of piecewise constant functions on some partition $(\Il)_{\lamm}$ of $\X$. 
The empirical risk minimizer $\ERM_m$ on $S_m$ is called a {\em regressogram}.
$S_m$ is a vector space of dimension $D_m = \card(\Lambda_m)$, spanned by the family $(\1_{\Il})_{\lamm}$. 
Since this basis is orthogonal in $L^2(\mu)$ for any probability measure  $\mu$ on $\X$, computations are quite easy. 
In particular, we have:
\[ \bayes_m = \sum_{\lamm} \betl \1_{\Il} \quad \mbox{and} \quad \ERM_m = \sum_{\lamm} \bethl \1_{\Il} \enspace ,\]
where 
\begin{gather*}
\betl \egaldef \E_P \croch{ Y \sachant X \in \Il } \qquad \bethl \egaldef \frac{1}{n\phl} \sum_{X_i \in \Il} Y_i \qquad \phl \egaldef P_n(X \in \Il) \enspace .
\end{gather*}
Note that $\ERM_m$ is uniquely defined if and only if each $\Il$ contains at least one of the $X_i$. 
Otherwise, $\ERM_m$ is not uniquely defined and we consider that the model $m$ cannot be chosen. 

\subsection{Main assumptions} \label{HP.sec.theorie.hyp}
In this section, we make the following assumptions. 
First, each model $S_m$ is a set of piecewise constants functions on some fixed partition $(\Il)_{\lamm}$ of $X$. 
Second, the family $(S_m)_{\mM_n}$ satisfies:
\begin{enumerate}
\item[\hypPpoly] Polynomial complexity of $\M_n$: $\card(\M_n) \leq \cM n^{\aM}$.
\item[\hypPrich] Richness of $\M_n$: $\exists m_0 \in \M_n$ s.t. $D_{m_0} \in \croch{\sqrt{n},c_{\mathrm{rich}} \sqrt{n}}$.
\end{enumerate}
Assumption \hypPpoly\ is quite classical for proving the asymptotic optimality of a model selection procedure; it is for instance implicitly assumed by \cite{KCLi:1987} in the homoscedastic fixed-design case.
Assumption \hypPrich\ is merely technical and can be changed if necessary; it only ensures that $(S_m)_{\mM_n}$ does not contain only models which are either too small or too large.

For any penalty function $\pen: \M_n \flens \R^+$, we define the following model selection procedure:
\begin{equation}
\mh \in \arg\min_{\mM_n, \, \min_{\lamm} \set{\phl} >0 } \set{ P_n \gamma(\ERM_m) + \pen(m) } \enspace . \label{HP.eq.mh.opt}
\end{equation}
Moreover, the data $(X_i,Y_i)_{1 \leq i \leq n}$ are assumed to be i.i.d. and to satisfy:
\begin{enumerate}
\item[\hypAb] The data are bounded: $\norm{Y_i}_{\infty} \leq A < \infty$.
\item[\hypAn] Uniform lower-bound on the noise level: $\sigma(X_i) \geq \sigmin>0$ a.s.
\item[\hypApu] The bias decreases as a power of $D_m$: there exist some $\betamaj,\cbiasmaj>0$ such that \[ \perte{\bayes_m} \leq \cbiasmaj D_m^{-\betamaj} \enspace . \]
\item[\hypArXl] Lower regularity of the partitions for $\loi(X)$: $ D_m \min_{\lamm} \set{\Prob\paren{X \in \Il}} \geq \crXl$.
\end{enumerate}
Further comments are made in Sections~\ref{HP.sec.theorie.mini} and~\ref{HP.sec.theorie.opt} about these assumptions, in particular about their possible weakening.

\subsection{Minimal penalties} \label{HP.sec.theorie.mini}
Our first result concerns the existence of a minimal penalty.
In this subsection, \hypPrich\ is replaced by the following strongest assumption:
\begin{enumerate}
\item[\hypPrichMin] $\exists c_0, c_{\mathrm{rich}} >0$ s.t. $\forall l \in \croch{\sqrt{n}, c_0 n/(c_{\mathrm{rich}} \ln(n))}$, 
$\exists m \in \M_n$ s.t. $D_{m} \in \croch{l,c_{\mathrm{rich}} l}$.
\end{enumerate}
The reason why \hypPrich\ is not sufficient to prove Theorem~\ref{HP.the.mini} below is that at least one model of dimension of order $n/\ln(n)$ should belong to the family $\paren{S_m}_{\mM_n}$; otherwise, it may not be possible to prove that such models are selected by penalization procedures beyond the minimal penalty.

\begin{theorem} \label{HP.the.mini}
Suppose all the assumptions of Section~\ref{HP.sec.theorie.hyp} are satisfied. 
Let $K \in [0;1)$, $L>0$, and assume that an event of probability at least $1 - L n^{-2}$ exists on which 
\begin{equation}  \label{HP.eq.pen.mini}
\forall \mM_n, \quad 0 \leq \pen(m) \leq K \E \croch{P_n \paren{ \gamma(\bayes_m) - \gamma(\ERM_m) }}  \enspace . 
\end{equation}

Then, there exist two positive constants $\Kminpr$, $\Kmindim$ such that, with probability at least $1 - \Kminpr n^{-2}$, 
\begin{equation*} 
D_{\mh} \geq \Kmindim n \ln(n)^{-1} \enspace ,
\end{equation*}
where $\mh$ is defined by \eqref{HP.eq.mh.opt}.
On the same event,
\begin{equation} \label{HP.eq.mini.risque}
\perte{\ERM_{\mh}} \geq \ln(n) \inf_{\mM_n} \set{\perte{\ERM_m}} \enspace .
\end{equation}

The constants $\Kminpr$ and $\Kmindim$ may depend on $K$, $L$ and constants in \hypPpoly, \hypPrichMin, \hypAb, \hypAn, \hypApu\ and \hypArXl, but do not depend on $n$. 
\end{theorem}
This theorem thus validates the first part of the heuristics of Section~\ref{HP.sec.cadre.heuristic}, proving that a minimal amount of penalization is required; when the penalty is smaller, the selected dimension $D_{\mh}$ and the quadratic risk of the final estimator $\perte{\ERM_{\mh}}$ blow up. 
This coupling is quite interesting, since the dimension $D_{\mh}$ is known in practice, contrary to $\perte{\ERM_{\mh}}$. It is then possible to detect from the data whether the penalty is too small, as proposed in Algorithm~\ref{HP.def.proc.gal.pente}.

The main interest of this result is its combination with Theorem~\ref{HP.the.opt} below. 
Nevertheless Theorem~\ref{HP.the.mini} is also interesting by itself for understanding the theoretical properties of penalization procedures. 
Indeed, it generalizes the results of \cite{Bir_Mas:2006} on the existence of minimal penalties to heteroscedastic regression with a random design, even if we have to restrict to regressograms. 
Moreover, we have a general formulation for the minimal penalty 
\[ \pen_{\min}(m) \egaldef \E \croch{P_n \paren{ \gamma(\bayes_m) - \gamma(\ERM_m) }} = \E\croch{p_2(m)} \enspace ,\]
which can be used in frameworks situations where it is not proportional to the dimension $D_m$ of the models (see Section~\ref{HP.sec.pratique.shape} and references therein). 

\medskip

In addition, assumptions \hypAb\ and \hypAn\ on the data are much weaker than the Gaussian homoscedastic assumption. 
They are also much more realistic, and moreover can be strongly relaxed. 
Roughly speaking, boundedness of data can be replaced by conditions on moments of the noise, and the uniform lower bound $\sigmin$ is no longer necessary when $\sigma$ satisfies some mild regularity assumptions. 
We refer to \citet[Section~4.3]{Arl:2008b} for detailed statements of these assumptions, and explanations on how to adapt proofs to these situations.

\medskip

Finally, let us comment on conditions \hypApu\ and \hypArXl. 
The upper bound \hypApu\ on the bias occurs in the most reasonable situations, for instance when $\X \subset \R^k$ is bounded, the partition $\paren{\Il}_{\lamm}$ is regular and the regression function $\bayes$ is $\alpha$-H\"olderian for some $\alpha>0$ ($\betamaj$ depending on $\alpha$ and $k$). 
It ensures that medium and large models have a significantly smaller bias than smaller ones; otherwise, the selected dimension would be allowed to be too small with significant probability.
On the other hand, \hypArXl\ is satisfied at least for ``almost regular'' partitions $\paren{\Il}_{\lamm}$, when $X$ has a lower bounded density w.r.t. the Lebesgue measure on $\X \subset \R^k$.

Theorem~\ref{HP.the.mini} is stated with a general formulation of \hypApu\ and \hypArXl, instead of assuming for instance that $\bayes$ is $\alpha$-H\"olderian and $X$ has a lower bounded density w.r.t $\Leb$, in order to point out the {\em generality} of the ``minimal penalization'' phenomenon. 
It occurs as soon as the models are not too much pathological. 
In particular, we do not make any assumption on the distribution of $X$ itself, but only that the models are not too badly chosen according to this distribution. 
Such a condition can be checked in practice if some prior knowledge on $\mathcal{L}(X)$ is available; if part of the data are unlabeled---a usual case---, classical density estimation procedures can be applied for estimating $\loi(X)$ from unlabeled data \citep{Dev_Lug:2001}. 

\subsection{Optimal penalties} \label{HP.sec.theorie.opt}
Algorithm~\ref{HP.def.proc.gal.pente} relies on a link between the minimal penalty pointed out by Theorem~\ref{HP.the.mini} and some optimal penalty. 
The following result is a formal proof of this link in the framework we consider: penalties close to twice the minimal penalty satisfy an oracle inequality with leading constant approximately equal to one.

\begin{theorem} \label{HP.the.opt}
Suppose all the assumptions of Section~\ref{HP.sec.theorie.hyp} are satisfied together with
\begin{enumerate}
\item[\hypAp] The bias decreases like a power of $D_m$: there exist $\betamin\geq\betamaj>0$ and $\cbiasmaj,\cbiasmin>0$ such that \[ \cbiasmin D_m^{-\betamin} \leq \perte{\bayes_m} \leq \cbiasmaj D_m^{-\betamaj} \enspace . \]
\end{enumerate}
Let $\delta \in (0,1)$, $L>0$, and assume that an event of probability at least $1 - L n^{-2}$ exists on which for every $\mM_n$,
\begin{equation}  \label{HP.eq.pen.opt}
(2-\delta) \E \croch{P_n \paren{ \gamma(\bayes_m) - \gamma(\ERM_m) }} \leq \pen(m) \leq (2+\delta) \E \croch{P_n \paren{ \gamma(\bayes_m) - \gamma(\ERM_m) }} \enspace . 
\end{equation}

Then, for every $0 < \eta < \min\set{\betamaj;1}/2$, there exist a constant $\Koptpr$ and a sequence $\epsilon_n$ tending to zero at infinity such that, with probability at least $1 - \Koptpr n^{-2}$, 
\begin{align} \notag 
D_{\mh} &\leq n^{1 - \eta} \\ \mbox{and} \qquad  
 \label{HP.eq.opt.traj}
\perte{\ERM_{\mh}} &\leq \paren{ \frac{1 + \delta}{1 - \delta} + \epsilon_n } \inf_{\mM_n} \set{ \perte{\ERM_m} } \enspace ,\end{align}
where $\mh$ is defined by \eqref{HP.eq.mh.opt}. 
Moreover, we have the oracle inequality
\begin{equation*} 
\E \croch{\perte{\ERM_{\mh}}} \leq \paren{ \frac{1 + \delta}{1 - \delta} + \epsilon_n } \E \croch{ \inf_{\mM_n} \set{ \perte{\ERM_m} } } + \frac{A^2 \Koptpr}{n^2} \enspace .
\end{equation*} 

The constant $\Koptpr$ may depend on $L,\delta,\eta$ and the constants in \hypPpoly, \hypPrich, \hypAb, \hypAn, \hypAp\ and \hypArXl, but not on $n$. 
The term $\epsilon_n$ is smaller than $\ln(n)^{-1/5}$; it can be made smaller than $n^{-\delta}$ for any $\delta \in (0;\delta_0(\betamin,\betamaj))$ at the price of enlarging $\Koptpr$.
\end{theorem}

This theorem shows that twice the minimal penalty $\pen_{\min}$ pointed out by Theorem~\ref{HP.the.mini} satisfies an oracle inequality with leading constant almost one. 
In other words, the slope heuristics of Section~\ref{HP.sec.cadre.heuristic} is valid.
The consequences of the {\em combination} of Theorems~\ref{HP.the.mini} and~\ref{HP.the.opt} are detailed in Section~\ref{HP.sec.theorie.comment}. 

The oracle inequality \eqref{HP.eq.opt.traj} remains valid when the penalty is only close to twice the minimal one.
In particular, {\em the shape of the penalty can be estimated by resampling} as suggested in Section~\ref{HP.sec.pratique.shape}. 

Actually, Theorem~\ref{HP.the.opt} above is a corollary of a more general result stated in Appendix~\ref{HP.sec.proofs.thm-gal}, Theorem~\ref{HP.the.opt.gal}. 
If 
\begin{equation} \label{HP.eq.pen.C} \pen(m) \approx K \E \croch{P_n \paren{ \gamma(\bayes_m) - \gamma(\ERM_m) }} \end{equation} instead of \eqref{HP.eq.pen.opt}, under the same assumptions, an oracle inequality with leading constant $C(K) + \epsilon_n$ instead of $1 + \epsilon_n$ holds with large probability. 
The constant $C(K)$ is equal to $(K-1)^{-1}$ when $K \in (1,2]$ and to $C(K)=K-1$ when $K>2$. 
Therefore, for every $K>1$, the penalty defined by \eqref{HP.eq.pen.C} is efficient up to a multiplicative constant.
This result is new in the heteroscedastic framework.

\medskip

Let us comment the additional assumption \hypAp, that is the lower bound on the bias. 
Assuming $\perte{\bayes_m}>0$ for every $\mM_n$ is classical for proving the asymptotic optimality of Mallows' $C_p$ \citep{Shi:1981,KCLi:1987,Bir_Mas:2006}. 
\hypAp\ has been made by \cite{Sto:1985} and \cite{Bur:2002} in the density estimation framework, for the same technical reasons as ours.
Assumption \hypAp\ is satisfied in several frameworks, such as the following: 
$(\Il)_{\lamm}$ is ``regular'', $X$ has a lower-bounded density w.r.t. the Lebesgue measure on $\X \subset \R^k$, and $\bayes$ is non-constant and $\alpha$-h\"olderian (w.r.t. $\norm{\cdot}_{\infty}$), with 
\[ \beta_1 = k^{-1} + \alpha^{-1} - (k-1) k^{-1} \alpha^{-1} \qquad \mbox{and} \qquad \beta_2 = 2 \alpha k^{-1} \enspace . \] 
We refer to \citet[Section~8.10]{Arl:2007:phd} for a complete proof. 

When the lower bound in \hypAp\ is no longer assumed, \eqref{HP.eq.opt.traj} holds with two modifications in its right-hand side \citep[for details, see][Remark~9]{Arl:2008b}: the $\inf$ is restricted to models of dimension larger than $\ln(n)^{\gamma_1}$, and there is a remainder term $\ln(n)^{\gamma_2} n^{-1}$, where $\gamma_1,\gamma_2>0$ are numerical constants. 
This is equivalent to \eqref{HP.eq.opt.traj}, unless there is a model of small dimension with a small bias.
The lower bound in \hypAp\ ensures that it cannot happen. 
Note that if there is a small model close to $\bayes$, it is hopeless to obtain an oracle inequality with a penalty which estimates $\penid$, simply because deviations of $\penid$ around its expectation would be much larger than the excess loss of the oracle. In such a situation, BIC-like methods are more appropriate; for instance, \cite{Csi:2002} and \cite{Csi_Shi:2000} showed that BIC penalties are minimal penalties for estimating the order of a Markov chain.

\subsection{Main theoretical and practical consequences} \label{HP.sec.theorie.comment}
The slope heuristics and the correctness of Algorithm~\ref{HP.def.proc.gal.pente} follow from the combination of Theorems~\ref{HP.the.mini} and~\ref{HP.the.opt}.

\subsubsection{Optimal and minimal penalties} \label{HP.sec.theorie.comment.pen}
For the sake of simplicity, let us consider the penalty $K \E\croch{p_2(m)}$ with any $K>0$; any penalty close to this one satisfies similar properties. 
At first reading, one can think of the homoscedastic case where $\E\croch{p_2(m)}\approx \sigma^2 D_m n^{-1}$; one of the novelties of our results is that the general picture is quite similar.

According to Theorem~\ref{HP.the.opt}, the penalization procedure associated with $K \E\croch{p_2(m)}$ satisfies an oracle inequality with leading constant $C_n(K)$ as soon as $K>1$, and $C_n(2) \approx 1$.
Moreover, results proved by \cite{Arl:2008a} imply that $C_n(K) \geq C(K) >1$ as soon as $K$ is not close to 2. 
Therefore, {\em $K=2$ is the optimal multiplying factor} in front of $\E\croch{p_2(m)}$.

When $K<1$, Theorem~\ref{HP.the.mini} shows that no oracle inequality can hold with leading constant $C_n(K) < \ln(n)$. 
Since $C_n(K) \leq (K-1)^{-1} < \ln(n)$ as soon as $K > 1 + \ln(n)^{-1}$, {\em $K=1$ is the minimal multiplying factor} in front of $\E\croch{p_2(m)}$.
More generally, $\pen_{\min}(m) \egaldef \E\croch{p_2(m)}$ is proved to be a {\em minimal penalty}.

\medskip

In short, Theorems~\ref{HP.the.mini} and~\ref{HP.the.opt} prove the slope heuristics described in Section~\ref{HP.sec.cadre.heuristic}:
\[ \mbox{\em ``optimal'' penalty} \approx 2 \times \mbox{\em ``minimal'' penalty} \enspace . \]
\cite{Bir_Mas:2006} have proved the validity of the slope heuristics in the Gaussian homoscedastic framework. 
This paper extends their result to a non-Gaussian and heteroscedastic setting.

\subsubsection{Dimension jump} \label{HP.sec.theorie.comment.dim}
In addition, Theorems~\ref{HP.the.mini} and~\ref{HP.the.opt} prove the existence of a crucial phenomenon: there exists a ``dimension jump''---complexity jump in the general framework---around the minimal penalty.
Let us consider again the penalty $K \E\croch{p_2(m)}$. 
As in Algorithm~\ref{HP.def.proc.gal.pente}, let us define
\[ \mh(K) \in \arg\min_{\mM_n} \set{ P_n \gamma\paren{\ERM_m} + K \E\croch{p_2(m)}} \enspace .\]
A careful look at the proofs of Theorems~\ref{HP.the.mini} and~\ref{HP.the.opt} shows that there exist constants $\Kdimjump,\Kdimjumpd >0$ and an event of probability $1 - \Kdimjump  n^{-2}$ on which
\begin{equation} \label{eq.dimjump} \forall 0<K <1-\frac{1}{\ln(n)}, \, D_{\mh(K)} \geq \frac{\Kdimjumpd n}{(\ln(n))^{2}}  
\,\, \mbox{ and } \, \,\forall K > 1+\frac{1}{\ln(n)}, \, D_{\mh(K)} \leq n^{1 - \eta} \enspace . 
\end{equation}
Therefore, the dimension $D_{\mh(K)}$ of the selected model jumps around the minimal value $K=1$, from values of order $n(\ln(n))^{-2}$ to $n^{1 - \eta}$. 

\medskip

Let us know explain why Algorithm~\ref{HP.def.proc.gal.pente} is correct, assuming that $\pen_{\mathrm{shape}}(m)$ is close to $\E\croch{p_2(m)}$. 
With definition \eqref{def.Kmin.thresh} of $\widehat{K}_{\min}$ and a threshold $D_{\mathrm{thresh}} \propto n (\ln(n))^{-3}$, \eqref{eq.dimjump} ensures that 
\[ 1 - \frac{1}{\ln(n)} \leq \widehat{K}_{\min} \leq 1 + \frac{1}{\ln(n)} \]
with a large probability. 
Then, according to Theorem~\ref{HP.the.opt}, the output of Algorithm~\ref{HP.def.proc.gal.pente} satisfies an oracle inequality with leading constant $C_n$ tending to one as $n$ tends to infinity.

\subsection{Comparison with data-splitting methods} \label{HP.sec.theorie.comment.cv}
Tuning parameters are often chosen by cross-validation or by another data-splitting method, which suffer from some drawbacks compared to Algorithm~\ref{HP.def.proc.gal.pente}.

First, $V$-fold cross-validation, leave-$p$-out and repeated learning-testing methods require a larger computation time. Indeed, they need to perform the empirical risk minimization process for each model several times, whereas Algorithm~\ref{HP.def.proc.gal.pente} only needs to perform it once.

Second, $V$-fold cross-validation is asymptotically suboptimal when $V$ is fixed, as shown by \cite[Theorem~1]{Arl:2008a}. 
The same suboptimality result is valid for the hold-out, when the size of the training set is not asymptotically equivalent to the sample size $n$. 
On the contrary, Theorems~\ref{HP.the.mini} and~\ref{HP.the.opt} prove that Algorithm~\ref{HP.def.proc.gal.pente} is asymptotically optimal in a framework including the one used by \cite[Theorem~1]{Arl:2008a} for proving the suboptimality of $V$-fold cross-validation. 
Hence, the quadratic risk of Algorithm~\ref{HP.def.proc.gal.pente} should be smaller, within a factor $\kappa>1$.

Third, hold-out with a training set of size $n_t \sim n$, for instance $n_t = n - \sqrt{n}$ or $n_t = n (1 - \ln(n)^{-1})$, is known to be unstable. 
The final output $\mh$ strongly depends on the choice of a particular split of the data. 
According to the simulation study of Section~\ref{HP.sec.algo.defKmin}, Algorithm~\ref{HP.def.proc.gal.pente} is far more stable.

\medskip

To conclude, compared to data splitting methods, Algorithm~\ref{HP.def.proc.gal.pente} is either faster to compute, more efficient in terms of quadratic risk, or more stable.
Then, Algorithm~\ref{HP.def.proc.gal.pente} should be preferred each time it can be used.
Another approach is to use aggregation techniques, instead of selecting one model. 
As shown by several results \citep[see for instance][]{Tsy:2004,Lec:2007a}, aggregating estimators built upon a training simple of size $n_t \sim n$ can have an optimal quadratic risk. 
Moreover, aggregation requires approximately the same computation time as Algorithm~\ref{HP.def.proc.gal.pente}, and is much more stable than the hold-out. 
Hence, it can be an alternative to model selection with Algorithm~\ref{HP.def.proc.gal.pente}.

\section{Conclusion} \label{HP.sec.conclu}
This paper provides mathematical evidence that the method introduced by \cite{Bir_Mas:2006} for designing data-driven penalties remains efficient in a non-Gaussian framework. 
The purpose of this conclusion is to relate the slope heuristics developed in Section~\ref{HP.sec.cadre} to the well known Mallows' $C_{p}$ and Akaike's criteria and to the unbiased estimation of the risk
principle. 

Let us come come back to Gaussian model selection in order to explain how to guess what is the right penalty from the data themselves.
Let $\gamma_{n}$ be some empirical criterion (for instance the least-squares criterion as in this paper, or the log-likelihood criterion), $\paren{  S_{m}}  _{\mM_n}$ be a collection of models and for every $\mM_n$ $\bayes_{m}$ be some minimizer  of $t\flapp \E\croch{ \gamma_{n} \paren{  t} }$ over $S_{m}$ (assuming that such a point exists).
Minimizing some penalized criterion
\[
\gamma_{n} \paren{\ERM_m } + \pen (m)
\]
over $\M_n$ amounts to minimize%
\begin{gather*}
\widehat{b}_{m} - \widehat{v}_{m} +\pen (m) \enspace ,
\\
\mbox{where} \quad \forall \mM_n, \,\, \widehat{b}_{m} = \gamma_{n} \paren{ \bayes_{m} }  - \gamma_{n} \paren{\bayes}
\mbox{ and }
\widehat{v}_{m} = \gamma_{n} \paren{ \bayes_{m} }  - \gamma_{n} \paren{ \ERM_m } \enspace .
\end{gather*}
The point is that $\widehat{b}_{m}$ is an unbiased estimator of the bias term $\perte{\bayes_m}$. 
Having concentration arguments in mind, minimizing $\widehat{b}_{m} - \widehat{v}_{m} +\pen (m)$ can be conjectured approximately equivalent to minimize
\[
\perte{\bayes_m} -\E\croch{  \widehat{v}_{m} } +\pen(m) \enspace .
\]
Since the purpose of model selection is to minimize the risk $\E\croch{\perte{\ERM_m}} $, an ideal penalty would
be
\[
\pen (m)  = \E\croch{ \widehat{v}_{m} } +\E\croch{ \pertedeux{\bayes_m}{\ERM_m} } \enspace .%
\]

In Gaussian least-squares regression with a fixed design, the models $S_{m}$ are linear and $\E\croch{  \widehat{v}_{m}}  =\E\croch{  \pertedeux{\bayes_m}{\ERM_m} }$ is explicitly computable if the noise level is constant and known; this leads to Mallows' $C_{p}$ penalty.
When $\gamma_n$ is the log-likelihood, 
\[
\E\croch{  \widehat{v}_{m} }  \approx \E\croch{  \pertedeux{\bayes_m}{\ERM_m} } \approx \frac{ D_{m}} {2n} 
\]
asymptotically, where $D_{m}$ stands for the number of parameters defining model $S_{m}$; this leads to Akaike's Information Criterion (AIC).
Therefore, both Mallows' $C_{p}$ and Akaike's criterion are based on the unbiased (or asymptotically unbiased) risk estimation principle.

This paper goes further in this direction, using that $\E\croch{  \widehat{v}_{m} }  \approx \E \croch{  \pertedeux{\bayes_m}{\ERM_m} }$ remains a valid approximation in a non-asymptotic framework. 
Then, a good penalty becomes $2\E \croch{  \widehat{v}_{m}}$ or $2\widehat{v}_{m}$, having in mind concentration arguments. 
Since $\widehat{v}_{m}$ is the minimal penalty, this explains the slope heuristics \citep{Bir_Mas:2006} and connects it to Mallows' $C_{p}$ and Akaike's heuristics.

\medskip

The second main idea developed in this paper is that the minimal penalty can be estimated from the data; Algorithm~\ref{HP.def.proc.gal.pente} uses the jump of complexity which occurs around the minimal penalty, as shown in Sections~\ref{HP.sec.algo.defKmin} and~\ref{HP.sec.theorie.comment.dim}.
Another way to estimate the minimal penalty when it is (at least approximately) of the form $\alpha D_{m}$ is to estimate $\alpha$ by the {\em slope} of the graph of $\gamma_{n}\paren{\ERM_m}$ for large enough values of $D_{m}$; this method can be extended to other shapes of penalties, simply by replacing $D_{m}$ by some (known!) function $f\paren{D_{m}}$. 

The slope heuristics can even be combined with resampling ideas, by taking a function $f$ built from a randomized empirical criterion. 
As shown by \cite{Arl:2008:shape}, this approach is much more efficient than the rougher choice $f\paren{D_{m}}=D_{m}$ for heteroscedastic regression frameworks. 
The question of the optimality of the slope heuristics in general remains an open problem; nevertheless, we believe that this heuristics can be useful in practice, and that proving its efficiency in this paper helps to understand it better.

\medskip

Let us finally mention that contrary to \cite{Bir_Mas:2006}, we assume in this paper that the collection of models $\M_n$ is ``small'', that is $\card(\M_n)$ grows at most like a power of $n$. 
For several problems, such that complete variable selection, larger collections of models have to be considered; then, it is known from the homoscedastic case that the minimal penalty is much larger than $\E\croch{p_2(m)}$.
Nevertheless, \'Emilie Lebarbier has used the slope heuristics with $f\paren{ D_{m}}  = D_{m} \paren{  2.5+\ln\paren{  \frac{n}{D_{m}}}}  $ for multiple change-points detection from $n$ noisy data, using the results by \cite{Bir_Mas:2006} in the Gaussian case.

Let us now explain how we expect to generalize the slope heuristics to the non-Gaussian heteroscedastic case when $\M_n$ is large.
First, group the models according to some complexity index $C_m$ such as their dimensions $D_m$; for $C \in \set{1, \ldots, n^k }$, define $\widetilde{S_C} = \bigcup_{C_m = C} S_m$. 
Then, replace the model selection problem with the family $(S_m)_{\mM_n}$ by a ``complexity selection problem'', that is model selection with the family $\paren{ \widetilde{S_C} }_{1 \leq C \leq n^k}$. 
We conjecture that this grouping of the models is sufficient to take into account the richness of $\M_n$ for the optimal calibration of the penalty. 
A theoretical justification of this point could rely on the extension of our results to any kind of model, since $\widetilde{S_C}$ is not a vector space in general. 


\acks{The authors gratefully acknowledge the anonymous referees for several suggestions and references.}

\appendix

\section{Proofs} \label{HP.sec.proofs}
This appendix is devoted to the proofs of the results stated in the paper. Proposition~\ref{HP.pro.algo.step} is proved in Section~\ref{HP.sec.algo.proof}; Theorem~\ref{HP.the.opt} is proved in Sections~\ref{HP.sec.proofs.thm-gal} and~\ref{sec.pr.the.opt.gal}; Theorem~\ref{HP.the.mini} is proved in Section~\ref{sec.pr.the.mini}; the remaining sections are devoted to probabilistic results used in the main proofs and technical proofs. 
\subsection{Conventions and notations} \label{HP.sec.proofs.conv}
In the rest of the paper, $L$ denotes a universal constant, not
necessarily the same at each occurrence.
When $L$ is not universal, but depends on $p_1, \ldots, p_k$, it is written $L_{p_1,\ldots, p_k}$.
Similarly, $L_{\hypHPm}$ (resp. $L_{\hypHP}$) denotes a constant allowed to depend on the parameters of the assumptions made in Theorem~\ref{HP.the.mini} (resp. Theorem~\ref{HP.the.opt.gal}), including \hypPpoly\ and \hypPrich. 
We also make use of the following notations: 
\begin{itemize}
\item $\forall a,b \in \R$, $\mini{a}{b}$ is the minimum of $a$ and $b$, $\maxi{a}{b}$ is the maximum of $a$ and $b$, $a_+ = \maxi{a}{0}$ is the positive part of $a$ and $a_- = \mini{a}{0}$ is its negative part.
\item $\forall \Il \subset \X$, $\pl \egaldef P(X \in \Il)$ and $\sigl^2 \egaldef \E\croch{\paren{Y - \bayes_m(X)}^2 \sachant X \in \Il}$.
\item Since $\E\croch{p_1(m)}$ is not well-defined (because of the event $\set{ \min_{\lamm} \set{\phl} = 0 }$), we have to take the following convention 
\begin{equation*}
p_1(m) = \punmin (m) \egaldef \sum_{\lamm \telque \phl>0} \pl \paren{\betl - \bethl}^2 + \sum_{\lamm \telque \phl=0} \pl \sigl^2 \enspace .
\end{equation*}
Remark that $p_1(m) = \punmin(m)$ when $\min_{\lamm} \set{\phl}>0$), so that this convention has no consequences on the final results (Theorems~\ref{HP.the.mini} and~\ref{HP.the.opt.gal}).
\end{itemize}

\subsection{Proof of Proposition~\ref{HP.pro.algo.step}} \label{HP.sec.algo.proof}
First, since $\M_n$ is finite, the infimum in \eqref{HP.def.algo.step.Ki} is attained as soon as $G(m_{i-1}) \neq \emptyset$, so that $m_i$ is well defined for every $i \leq i_{\max}$.
Moreover, by construction, $g(m_i)$ decreases with $i$, so that all the $m_i \in \M_n$ are different; hence, Algorithm~\ref{HP.algo.step} terminates and $i_{\max} +1 \leq \card(\M_n)$.
We now prove by induction the following property for every $i \in \set{0, \ldots, i_{\max}  }$:
\[ \mathcal{P}_i : \qquad K_i < K_{i+1} \quad \mbox{ and } \quad \forall K \in [K_i, K_{i+1}), \quad \mh(K) = m_i \enspace . \]
Notice also that $K_i$ can always be defined by \eqref{HP.def.algo.step.Ki} with the convention $\inf\emptyset=+\infty$.

\medskip
\subsubsection*{$\mathcal{P}_0$ holds true} 
By definition of $K_1$, it is clear that $K_1>0$ (it may be equal to $+\infty$ if $G(m_0)=\emptyset$).
For $K=K_0=0$, the definition of $m_0$ is the one of $\mh(0)$, so that $\mh(K) = m_0$.
For $K \in (0,K_1)$, Lemma~\ref{HP.le.algo.step} shows that either $\mh(K)=\mh(0)=m_0$ or $\mh(K) \in G(0)$. In the latter case, by definition of $K_1$, 
\[ \frac{f(\mh(K)) - f(m_0)} {g(m_0) - g(\mh(K))} \geq K_1 > K \]
hence
\[ f(\mh(K)) + K g(\mh(K)) > f(m_0) + K g(m_0) \] which is contradictory with the definition of $\mh(K)$. 
Therefore, $\mathcal{P}_0$ holds true.

\medskip
\subsubsection*{$\mathcal{P}_i \Rightarrow \mathcal{P}_{i+1}$ for every $i \in \set{0, \ldots, i_{\max}-1}$}
Assume that $\mathcal{P}_i$ holds true. 
First, we have to prove that $K_{i+2}>K_{i+1}$. Since $K_{i_{\max} + 1}=+\infty$, this is clear if $i=i_{\max}-1$. Otherwise, $K_{i+2}<+\infty$ and $m_{i+2}$ exists. Then, by definition of $m_{i+2}$ and $K_{i+2}$ (resp. $m_{i+1}$ and $K_{i+1}$), we have
\begin{align}
\label{HP.eq.pro.algo.step.1}
f(m_{i+2}) - f(m_{i+1}) = K_{i+2} ( g(m_{i+1}) - g(m_{i+2}))
\\
\label{HP.eq.pro.algo.step.2}
f(m_{i+1}) - f(m_{i}) = K_{i+1} ( g(m_{i}) - g(m_{i+1})) \enspace .
\end{align}
Moreover, $m_{i+2} \in G(m_{i+1}) \subset G(m_i)$, and $m_{i+2} \prec m_{i+1}$ (because $g$ is non-decreasing). Using again the definition of $K_{i+1}$, we have\begin{equation}
\label{HP.eq.pro.algo.step.3}
f(m_{i+2}) - f(m_i) > K_{i+1} (g(m_i) - g(m_{i+2}))
\end{equation}
(otherwise, we would have $m_{i+2} \in F_{i+1}$ and $m_{i+2} \prec m_{i+1}$, which is not possible).
Combining the difference of \eqref{HP.eq.pro.algo.step.3} and \eqref{HP.eq.pro.algo.step.2} with \eqref{HP.eq.pro.algo.step.1}, we have
\[ K_{i+2} (g(m_{i+1}) - g(m_{i+2})) > K_{i+1} (g(m_{i+1}) - g(m_{i+2})) \enspace , \]
hence $K_{i+2} > K_{i+1}$, since $g(m_{i+1}) > g(m_{i+2})$.

Second, we prove that $\mh(K_{i+1})=m_{i+1}$. From $\mathcal{P}_i$, we know that for every $\mM_n$, for every $K\in[K_i, K_{i+1})$, $f(m_i) + K g(m_i) \leq f(m) + K g(m)$. Taking the limit when $K$ tends to $K_{i+1}$, it follows that $m_i \in E(K_{i+1})$. By \eqref{HP.eq.pro.algo.step.2}, we then have $m_{i+1} \in E(K_{i+1})$. On the other hand, if $m \in E(K_{i+1})$, Lemma~\ref{HP.le.algo.step} shows that either $f(m) = f(m_i)$ and $g(m) = g(m_i)$ or $m \in G(m_i)$. In the first case, $m_{i+1} \prec m$ (because $g$ is non-decreasing). In the second one, $m \in F_{i+1}$, so $m_{i+1} \preceq m$. Since $\mh(K_{i+1})$ is the smallest element of $E(K_{i+1})$, we have proved that $m_{i+1} = \mh(K_{i+1})$.

Last, we have to prove that $\mh(K)=m_{i+1}$ for every $K \in (K_1,K_2)$. From the last statement of Lemma~\ref{HP.le.algo.step}, we have either $\mh(K)=\mh(K_1)$ or $\mh(K_1) \in G(\mh(K))$. In the latter case (which is only possible if $K_{i+2} < \infty$), by definition of $K_{i+2}$, 
\[ \frac{f(\mh(K)) - f(m_{i+1})} {g(m_{i+1}) - g(\mh(K))} \geq K_{i+2} > K \]
so that 
\[ f(\mh(K)) + K g(\mh(K)) > f(m_{i+1}) + K g(m_{i+1}) \] which is contradictory with the definition of $\mh(K)$. \BlackBox

\begin{lemma} \label{HP.le.algo.step}
With the notations of Proposition~\ref{HP.pro.algo.step} and its proof, if $0  \leq K < K^{\prime}$, $m \in E(K)$ and $m^{\prime} \in E(K^{\prime})$, then one of the two following statements holds true:
\begin{enumerate}
\item[(a)] $f(m) = f(m^{\prime})$ and $g(m) = g(m^{\prime})$.
\item[(b)] $f(m) < f(m^{\prime})$ and $g(m)>g(m^{\prime})$.
\end{enumerate}
In particular, either $\mh(K)=\mh(K^{\prime})$ or $\mh(K^{\prime}) \in G(\mh(K))$.
\end{lemma}
\begin{proof}
By definition of $E(K)$ and $E(K^{\prime})$, 
\begin{align}
\label{HP.eq.le.algo.step.1}
f(m) + K g(m) &\leq f(m^{\prime}) + K g(m^{\prime}) \\
f(m^{\prime}) + K^{\prime} g(m^{\prime}) &\leq f(m) + K^{\prime} g(m) \enspace .
\label{HP.eq.le.algo.step.2}
\end{align}
Summing \eqref{HP.eq.le.algo.step.1} and \eqref{HP.eq.le.algo.step.2} gives $(K^{\prime} - K) g(m^{\prime}) \leq (K^{\prime} - K) g(m)$ so that 
\begin{equation}
\label{HP.eq.le.algo.step.3} g(m^{\prime}) \leq g(m) \enspace . \end{equation}
Since $K \geq 0$, \eqref{HP.eq.le.algo.step.1} and \eqref{HP.eq.le.algo.step.3} give $f(m) + K g(m) \leq f(m^{\prime}) + K g(m)$, that is
\begin{equation}
\label{HP.eq.le.algo.step.4} f(m) \leq f(m^{\prime}) \enspace . \end{equation}

Moreover, \eqref{HP.eq.le.algo.step.4} and \eqref{HP.eq.le.algo.step.2} imply $g(m) = g(m^{\prime})$, hence  $f(m^{\prime}) \leq f(m)$, that is $f(m) = f(m^{\prime})$ by \eqref{HP.eq.le.algo.step.4}. 
Similarly, \eqref{HP.eq.le.algo.step.1} and \eqref{HP.eq.le.algo.step.3} show that $f(m) = f(m^{\prime})$ imply $g(m) = g(m^{\prime})$. In both cases, (a) is satisfied.
Otherwise, $f(m) < f(m^{\prime})$ and $g(m)>g(m^{\prime})$, that is the (b) statement.

The last statement follows by taking $m=\mh(K)$ and $m^{\prime}=\mh(K^{\prime})$, because $g$ is non-decreasing, so that the minimum of $g$ in $E(K)$ is attained by $\mh(K)$.
\end{proof}

\subsection{A general oracle inequality} \label{HP.sec.proofs.thm-gal}
First of all, let us state a general theorem, from which Theorem~\ref{HP.the.opt} is an obvious corollary.

\begin{theorem} \label{HP.the.opt.gal}
Suppose all the assumptions of Section~\ref{HP.sec.theorie.hyp} are satisfied together with
\begin{enumerate}
\item[\hypAp] The bias decreases like a power of $D_m$: there exist $\betamin\geq\betamaj>0$ and $\cbiasmaj,\cbiasmin>0$ such that \[ \cbiasmin D_m^{-\betamin} \leq \perte{\bayes_m} \leq \cbiasmaj D_m^{-\betamaj} \enspace . \]
\end{enumerate}

Let $L,\xi,c_1,C_1,C_2 \geq 0$, $c_2 > 1$ and assume that an event of probability at least $1 - L n^{-2}$ exists on which, for every $\mM_n$ such that $D_m \geq \ln(n)^{\xi}$,
\begin{equation}  \label{HP.eq.pen.opt.gal}
\begin{split}
& \E \croch{c_1 P \paren{ \gamma(\ERM_m) - \gamma(\bayes_m) } + c_2 P_n \paren{ \gamma(\bayes_m) - \gamma(\ERM_m) }}  \\
& \quad \leq \pen(m) \leq \E \croch{C_1 P \paren{ \gamma(\ERM_m) - \gamma(\bayes_m) } + C_2 P_n \paren{ \gamma(\bayes_m) - \gamma(\ERM_m) }} \enspace . \end{split} \end{equation}

Then, for every $0 < \eta < \min\set{\betamaj;1}/2$, there exist a constant $\Koptpr$ and a sequence $\epsilon_n$ tending to zero at infinity such that, with probability at least $1 - \Koptpr n^{-2}$,  
\begin{align} 
\notag D_{\mh} &\leq n^{1 - \eta}
\\
\label{HP.eq.opt.traj.gal}
\mbox{and} \qquad \perte{\ERM_{\mh}} &\leq \croch{ \frac{1+ (C_1 + C_2 - 2)_+}{\mini{\paren{c_1 + c_2 - 1}}{1}} + \epsilon_n } \inf_{\mM_n} \set{ \perte{\ERM_m} } 
\end{align}
where $\mh$ is defined by \eqref{HP.eq.mh.opt}.
Moreover, we have the oracle inequality
\begin{equation} \label{HP.eq.opt.class.gal}
\E \croch{\perte{\ERM_{\mh}}} \leq \croch{ \frac{1+ (C_1 + C_2 - 2)_+}{\mini{\paren{c_1 + c_2 - 1}}{1}} + \epsilon_n } \E \croch{ \inf_{\mM_n} \set{ \perte{\ERM_m} } } + \frac{A^2 \Koptpr}{n^2} \enspace .
\end{equation} 

The constant $\Koptpr$ may depend on $L$, $\eta$, $\xi$, $c_1$, $c_2$, $C_1$, $C_2$ and constants in \hypPpoly, \hypPrich, \hypAb, \hypAn, \hypAp\ and \hypArXl, but not on $n$. The term $\epsilon_n$ is smaller than $\ln(n)^{-1/5}$; it can be made smaller than $n^{-\delta}$ for any $\delta \in (0;\delta_0(\betamin,\betamaj))$ at the price of enlarging $\Koptpr$.
\end{theorem}

The particular form of condition \eqref{HP.eq.pen.opt.gal} on the penalty is motivated by the fact that the ideal shape of penalty $\E\croch{\penid(m)}$ (or equivalently $\E\croch{2 p_2(m)}$) is unknown in general. 
Then, it has to be estimated from the data, for instance by resampling. 
Under the assumptions of Theorem~\ref{HP.the.opt.gal}, \cite{Arl:2008a,Arl:2008b} has proved that resampling and $V$-fold penalties satisfy condition \eqref{HP.eq.pen.opt.gal} with constants $c_1 + c_2 = 2 - \delta_n$, $C_1 + C_2 = 2 + \delta_n$ (for some absolute sequence $\delta_n$ tending to zero at infinity), and some numerical constant $\xi>0$.
Then, Theorem~\ref{HP.the.opt.gal} shows that such a penalization procedure satisfies an oracle inequality with leading constant tending to 1 asymptotically.

\medskip

The rationale behind Theorem~\ref{HP.the.opt.gal} is that if $\pen(m)$ is close to $c_1 p_1(m) + c_2 p_2(m)$, then $\crit(m) \approx \perte{\bayes_m} + c_1 p_1(m) + (c_2 - 1) p_2(m)$. 
When $c_1 = c_2 = 1$, this is exactly the ideal criterion $\perte{\ERM_m}$. 
When $c_1 + c_2 = 2$ with $c_1 \geq 0$ and $c_2 > 1$, we obtain the same result because $p_1(m)$ and $p_2(m)$ are quite close, at least when $D_m$ is large enough. 
The closeness between $p_1$ and $p_2$ is the keystone of the slope heuristics. 
Notice that if $\max_{\mM_n} D_m \leq \Koptpr^{\prime} (\ln(n))^{-1} n$ (for some constant $\Koptpr^{\prime}$ depending only on the assumptions of Theorem~\ref{HP.the.opt}, as $\Koptpr$), one can replace the condition $c_2 > 1$ by $c_1 + c_2 > 1$ and $c_1,c_2\geq 0$ .

\subsection{Proof of Theorem~\ref{HP.the.opt.gal}} \label{sec.pr.the.opt.gal}
This proof is similar to the one of \citet[Theorem~1]{Arl:2008b}. We give it for the sake of completeness.

From \eqref{HP.eq.etape1}, we have for each $\mM_n$ such that $A_n(m) \egaldef \min_{\lamm} \set{n\phl}>0$
\begin{equation} 
\label{HP.eq.oracle.pr.1}
\perte{\ERM_{\mh}} - \paren{\penid^{\prime}(\mh) - \pen(\mh)} \leq \perte{\ERM_m} + \paren{\pen(m) - \penid^{\prime}(m)} \enspace . \end{equation}
with $\penid^{\prime}(m) \egaldef p_1(m) + p_2(m) - \delc(m) = \pen(m) + (P - P_n) \gamma(\bayes)$ and $\delc(m) \egaldef (P_n - P)(\gamma\paren{\bayes_m} - \gamma\paren{\bayes})$. It is sufficient to control $\pen - \penid^{\prime}$ for every $\mM_n$. 

We will thus use the concentration inequalities of Section~\ref{HP.sec.conc} with $x=\gamma \ln(n)$ and $\gamma=2+\aM$. Define $B_n(m) = \min_{\lamm} \set{n\pl}$, and $\Omega_{n}$ the event on which 
\begin{itemize}
\item for every $\mM_n$, \eqref{HP.eq.pen.opt.gal} holds 
\item for every $\mM_n$ such that $B_n(m) \geq 1$, \eqref{HP.eq.conc.p1.min} and \eqref{HP.eq.conc.p1.maj} hold:
\begin{align*} 
\punmin(m) &\geq \E \croch{\punmin(m)} - L_{\hypHP}  \croch{  \frac{\ln(n)^{2}}{\sqrt{D_m}} +  e^{-L B_n(m)} }  \E\croch{p_2(m)}
\\
\punmin(m) &\leq \E \croch{\punmin(m)} + L_{\hypHP}  \croch{  \frac{\ln(n)^{2}}{\sqrt{D_m}}  +  \sqrt{D_m} e^{-L B_n(m)} }  \E\croch{p_2(m)} 
\end{align*}
\item for every $\mM_n$ such that $B_n(m)>0$, \eqref{HP.eq.conc.p1.min.2}, \eqref{HP.eq.conc.p2.Qmp} and~\ref{HP.eq.conc.delta.borne.2} hold:
\begin{gather*} 
\punmin(m) \geq \paren{\frac{1}{2 + (\gamma+1) B_n(m)^{-1} \ln(n) } -  \frac{L_{\hypHP} 
\ln(n)^{2}}{\sqrt{D_m}} } \E\croch{p_2(m)}  \\
\absj{p_2(m) - \E\croch{p_2(m)}} \leq \frac{L_{\hypHP} \ln(n)}{\sqrt{D_m}} \croch{\perte{\bayes_m} + \E\croch{p_2(m)} } \\ 
 \absj{\delc(m)} \leq \frac{\perte{\bayes_m}}{\sqrt{D_m}} + L_{\hypHP}  \frac{\ln(n)}{\sqrt{D_m} }  \E\croch{p_2(m)}
\end{gather*}
\end{itemize}
From Proposition~\ref{HP.pro.conc.p1} (for $\punmin$),
Proposition~\ref{HP.pro.conc.p2} (for $p_2$) and Proposition~\ref{HP.pro.conc.delta.borne} (for $\delc(m)$), 
\[ \Prob\paren{\Omega_{n}} \geq 1 - L \sum_{\mM_n} n^{-2 - \aM} \geq 1 - L_{\cM} n^{-2} \enspace . \]

\medskip

For every $\mM_n$ such that $D_m \leq L_{\crXl} n \ln(n)^{-1}$, \hypArXl\ implies that $B_n(m) \geq L^{-1}\ln(n) \geq 1$.
As a consequence, on $\Omega_{n}$, if $\ln(n)^{7} \leq D_m \leq L_{\crXl} n \ln(n)^{-1}$:
\begin{equation*}
\begin{split}
\max \set{ \absj{\punmin(m) - \E\croch{\punmin(m)}} , \absj{p_2(m) - \E\croch{p_2(m)}} , \absj{\delc(m)}}  
 \leq \frac{L_{\hypHP} \E\croch{\perte{\bayes_m} + p_2(m)}} {\ln(n)}
\end{split}
\end{equation*}
Using \eqref{HP.eq.comp.Ep1.Ep2} (in Proposition~\ref{HP.pro.esp.p1p2}) and the fact that $B_n(m) \geq L^{-1} \ln(n)$, 
\[ \frac{(c_1 + c_2) \paren{1 - \widetilde{\delta_n}}}{2}  \leq \E\croch{\pen(m)} \leq \frac{(C_1 + C_2) \paren{1 + \widetilde{\delta_n}}} {2} \E\croch{\punmin(m)+p_2(m)} \] with $0 \leq \widetilde{\delta_n} \leq L \ln(n)^{-1/4}$.
We deduce: if $n \geq L_{\hypHP}$, for every $\mM_n$ such that $\ln(n)^{7} \leq D_m \leq L_{\crXl} n \ln(n)^{-1}$, on $\Omega_{n}$, 
\begin{equation*}
\begin{split}
 \croch{ \paren{c_1 + c_2 - 2}_- - \frac{L_{\hypHP}}{\ln(n)^{1/4}} } p_1(m) \leq (\pen-\penid^{\prime})(m) \\
 \leq \croch{ \paren{C_1 + C_2 - 2}_+  + \frac{L_{\hypHP}}{\ln(n)^{1/4}} } p_1(m) \enspace .\end{split}
\end{equation*}
We need to assume that $n$ is large enough in order to upper bound $\E\croch{p_2(m)}$ in terms of $p_1(m)$, since we only have \[ p_1(m) \geq \croch{ 1 - \frac{L_{\hypHP}}{\ln(n)^{1/4}} }_+ \E\croch{p_2(m)} \] in general.
Combined with \eqref{HP.eq.oracle.pr.1}, this gives: if $n \geq L_{\hypHP}$,
\begin{equation*} 
\begin{split}
\perte{\ERM_{\mh}} \1_{\ln(n)^{5} \leq D_{\mh} \leq L_{\crXl} n \ln(n)^{-1}} \leq \croch{ \frac{1+ (C_1 + C_2 - 2)_+}{\mini{\paren{c_1 + c_2 - 1}}{1}} + \frac{L_{\hypHP}}{\ln(n)^{1/4}} } \\
\times \inf_{\mM_n \telque \ln(n)^{7} \leq D_m \leq L_{\aM,\crXl} n \ln(n)^{-1}} \set{ \perte{\ERM_m} } \enspace .
\end{split}
\end{equation*}

We now use Lemmas~\ref{HP.le.borne.Dmh} and~\ref{HP.le.borne.Dmo} below to control on $\Omega_n$ the dimensions of the selected model $\mh$ and the oracle model $\mo \in \arg\min_{\mM_n} \set{\perte{\ERM_m}}$.

The result follows since $L_{\hypHP} \ln(n)^{-1/4} \leq \epsilon_n = \ln(n)^{-1/5}$ for $n \geq L_{\hypHP}$. We finally remove the condition $n \geq n_0 = L_{\hypHP}$ by choosing $\Koptpr = L_{\hypHP}$ such that $\Koptpr n_0^{-2} \geq 1$.

\paragraph{Classical oracle inequality}
Since \eqref{HP.eq.opt.traj.gal} holds true on $\Omega_n$, 
\begin{align*}
\E \croch{ \perte{\ERM_{\mh}} } &= \E \croch{ \perte{\ERM_{\mh}} \1_{\Omega_n} } + \E \croch{ \perte{\ERM_{\mh}} \1_{\Omega_n^c} } \\
&\leq \croch{ 2\eta - 1 + \epsilon_n }  \E \croch{ \inf_{\mM_n} \set{ \perte{\ERM_m} } } + A^2 \Koptpr \Prob\paren{\Omega_n^c} \end{align*}
which proves \eqref{HP.eq.opt.class.gal}. \BlackBox

\begin{lemma}[Control on the dimension of the selected model] \label{HP.le.borne.Dmh}
Let $c>0$ and $\alpha>\paren{1 - \betamaj}_+/2$. 
Then, if $n \geq L_{\hypHP,c,\alpha}$, on the event $\Omega_n$ defined in the proof of Theorem~\ref{HP.the.opt.gal}, 
\begin{equation*} 
\ln(n)^{7} \leq D_{\mh} \leq n^{1/2 + \alpha} \leq c n \ln(n)^{-1} \enspace .
\end{equation*}
\end{lemma}

\begin{lemma}[Control on the dimension of the oracle model] \label{HP.le.borne.Dmo}
Define the oracle model $\mo \in \arg\min_{\mM_n} \set{\perte{\ERM_m}}$.
Let $c>0$ and $\alpha>\paren{1 - \betamaj}_+/2$. 
Then, if $n \geq L_{\hypHP,c,\alpha}$, on the event $\Omega_n$ defined in the proof of Theorem~\ref{HP.the.opt.gal}, 
\begin{equation*} 
\ln(n)^{7} \leq D_{\mo} \leq n^{1/2 + \alpha} \leq c n \ln(n)^{-1}  \enspace .
\end{equation*}
\end{lemma}

\paragraph{Proof of Lemma~\ref{HP.le.borne.Dmh}}
By definition, $\mh$ minimizes $\crit(m)$ over $\M_n$. It thus also minimizes 
\[ \crit^{\prime}(m) = \crit(m) - P_n \gamma(\bayes) = \perte{\bayes_m} - p_2(m) + \delc(m) + \pen(m) \] over $\M_n$.
\begin{enumerate}
\item Lower bound on $\crit^{\prime}(m)$ for small models: let $\mM_n$ such that $D_m < \paren{\ln(n)}^{7}$. We then have
\begin{align*}
\perte{\bayes_m} &\geq \cbiasmin \paren{\ln(n)}^{- 7 \betamin} \qquad \mbox{from \hypAp} \\
\pen(m) &\geq 0  \\
p_2(m) &\leq L_{\hypHP} \sqrt{\frac{\ln(n)}{n} } + L_{\hypHP} \frac{D_m}{n} \leq L_{\hypHP} \sqrt{\frac{\ln(n)}{n} } \qquad \mbox{from \eqref{HP.eq.conc.p2}} \end{align*}
and from \eqref{HP.eq.conc.delta.borne.2} (in Proposition~\ref{HP.pro.conc.delta.borne}), 
\[ \delc(m) \geq - L_A \sqrt{\frac{\perte{\bayes_m} \ln(n)}{n}} + L_A \frac{\ln(n)}{n} \geq - L_A \sqrt{\frac{\ln(n)}{n}} \enspace .\] 
We then have \[ \crit^{\prime}(m) \geq L_{\hypHP} \paren{\ln(n)}^{-L_{\betamin}} \enspace . \]
\item Lower bound for large models: let $\mM_n$ such that $D_m \geq n^{1/2 + \alpha}$. From \eqref{HP.eq.pen.opt.gal} and \eqref{HP.eq.conc.p2} (in Proposition~\ref{HP.pro.conc.p2}),
\begin{align*}
\pen(m)-p_2(m) &\geq \paren{ c_2 - 1 } \E\croch{p_2(m)} - L_A \sqrt{ \frac{\ln(n)} {n} } \\
&\geq \frac{ (c_2 - 1 ) \sigmin^2 D_m}{n}  - L_A \sqrt{ \frac{\ln(n)} {n} } \end{align*}
and from \eqref{HP.eq.conc.delta.borne},
\begin{equation*}
\delc(m) \geq - L_{\hypHP} \sqrt{\frac{\ln(n)}{n}} \enspace .
\end{equation*}
Hence, if $D_m \geq n^{1/2 + \alpha}$ and $n \geq L_{\hypHP,\alpha}$
\begin{align*}
\crit^{\prime}(m) \geq \pen(m) + \delc(m) - p_2(m) \geq L_{\hypHP,\alpha} n^{-1/2 + \alpha} \enspace .
\end{align*}
\item There exists a better model for $\crit(m)$: from \hypPrich, there exists $m_0 \in \M_n$ such that $\sqrt{n} \leq D_{m_0} \leq c_{\mathrm{rich}} \sqrt{n}$. If moreover $n \geq L_{c_{\mathrm{rich}}, \alpha}$, then
\[ \ln(n)^{7} \leq \sqrt{n} \leq D_{m_0} \leq c_{\mathrm{rich}} \sqrt{n} \leq n^{1/2 + \alpha} \enspace .\]
By \eqref{HP.eq.threshold} in Lemma~\ref{HP.le.threshold}, $A_n(m_0) \geq 1$ with probability at least $1 - L n^{-2}$. \\
Using \hypAp, \[ \perte{\bayes_{m_0}} \leq \cbiasmaj c_{\mathrm{rich}}^{\betamaj} n^{- \betamaj / 2}   \] 
so that, when $n \geq L_{\hypHP}$,
\begin{align*}
\crit^{\prime}(m_0) &\leq \perte{\bayes_{m_0}} + \absj{\delc(m)} + \pen(m) \\
&\leq L_{\hypHP} \paren{n^{-\betamaj/2} + n^{-1/2}} \enspace .\end{align*}
If $n \geq L_{\hypHP,\alpha}$, this upper bound is smaller than the previous lower bounds for small and large models. \BlackBox
\end{enumerate}

\paragraph{Proof of Lemma~\ref{HP.le.borne.Dmo}}
Recall that $\mo$ minimizes $\perte{\ERM_m}=\perte{\bayes_m} + p_1(m)$ over $\mM_n$, with the convention $\perte{\ERM_m} = \infty$ if $A_n(m)=0$.
\begin{enumerate}
\item Lower bound on $\perte{\ERM_m}$ for small models: let $\mM_n$ such that $D_m < \paren{\ln(n)}^{7}$. From \hypAp, we have
\[ \perte{\ERM_m} \geq \perte{\bayes_m} \geq \cbiasmin \paren{\ln(n)}^{-7 \betamin} \enspace . \]
\item Lower bound on $\perte{\ERM_m}$ for large models: let $\mM_n$ such that $D_m > n^{1/2 + \alpha}$.
From \eqref{HP.eq.conc.p1.min.2}, for $n \geq L_{\hypHP,\alpha}$,
\begin{gather*}
\punmin(m) \geq  \paren{\frac{1}{2 + (\gamma+1) \paren{\crXl}^{-1} \ln(n) } -  \frac{L_{\hypHP,\alpha}} {n^{1/4}}  } \E\croch{\widetilde{p_2}(m)}  \\
\mbox{so that} \quad \perte{\ERM_m} \geq \punmin(m) \geq   L_{\hypHP,\alpha} n^{-1/2 + \alpha} \enspace .
\end{gather*}
\item There exists a better model for $\perte{\ERM_m}$: let $m_0 \in \M_n$ be as in the proof of Lemma~\ref{HP.le.borne.Dmh} and assume that $n \geq L_{c_{\mathrm{rich}},\alpha}$. Then,
\[ p_1(m_0) \leq L_{\hypHP} \E\croch{p_2(m)} \leq  L_{\hypHP} n^{-1/2} \]
and the arguments of the previous proof show that
\[\perte{\ERM_{m_0}} \leq L_{\hypHP} \paren{n^{- \betamaj / 2} + n^{-1/2}} \] which is smaller than the previous upper bounds for $n \geq L_{\hypHP,\alpha}$. \BlackBox
\end{enumerate}

\subsection{Proof of Theorem~\ref{HP.the.mini}} \label{sec.pr.the.mini}
Similarly to the proof of Theorem~\ref{HP.the.opt.gal}, we consider the event $\Omega^{\prime}_n$, of probability at least $1 - L_{\cM} n^{-2}$, on which:
\begin{itemize}
\item for every $\mM_n$, \eqref{HP.eq.pen.mini} (for $\pen$), \eqref{HP.eq.conc.p1.min.2} (for $\punmin$), \eqref{HP.eq.conc.p2}--\eqref{HP.eq.conc.p2.Qmp} (for $p_2$, with $x=\gamma\ln(n)$ and $\theta = \sqrt{\ln(n)/n}$) and \eqref{HP.eq.conc.delta.borne}--\eqref{HP.eq.conc.delta.borne.2} (for $\delc$, with $x=\gamma\ln(n)$ and $\eta = \sqrt{\ln(n)/n}$) hold true.
\item for every $\mM_n$ such that $B_n(m) \geq 1$, 
\eqref{HP.eq.conc.p1.min} and \eqref{HP.eq.conc.p1.maj} hold (for $\punmin$).
\end{itemize}

\paragraph{Lower bound on $D_{\mh}$}
By definition, $\mh$ minimizes \[ \crit^{\prime}(m) = \crit(m) - P_n \gamma(\bayes) = \perte{\bayes_m} - p_2(m) + \delc(m) + \pen(m) \] over $\mM_n$ such that $A_n(m) \geq 1$.
As in the proof of Theorem~\ref{HP.the.opt.gal}, we define $c = L_{\crXl} >0$ such that for every model of dimension $D_m \leq c n \ln(n)^{-1}$, $B_n(m) \geq L^{-1} \ln(n) \geq 1$. Let $c^{\prime} = \min(c,c_0)$ and $d \in (0,1)$ a constant to be chosen later.

\begin{enumerate}
\item Lower bound on $\crit^{\prime}(m)$ for ``small'' models: assume that $\mM_n$ and $D_m \leq d c^{\prime} n \ln(n)^{- 1}$. Then, $\perte{\bayes_m} + \pen(m) \geq 0$ and 
from \eqref{HP.eq.conc.delta.borne}, 
\[ \delc(m) \geq - L_A \sqrt{\frac{\ln(n)}{n}} \enspace .\] 
If $D_m \geq \ln(n)^4$, \eqref{HP.eq.conc.p2.Qmp} implies that 
\[ p_2(m) \leq \paren{1 + \frac{L_{\hypHPm}}{\ln(n)} } \E\croch{p_2(m)} \leq \frac{ L_{\hypHPm} D_m}{n} \leq \frac{c^{\prime} d  L_{\hypHPm}}{\ln(n)}  \enspace .\]
On the other hand, if $D_m < \ln(n)^{4}$, \eqref{HP.eq.conc.p2} implies that 
\[ p_2(m) \leq L_{\hypHPm} \sqrt{ \frac{\ln(n)}{n}} \enspace .\]
We then have \[ \crit^{\prime}(m) \geq - d L_{\hypHPm} \paren{\ln(n)}^{-1} \enspace . \]
\item There exists a better model for $\crit(m)$: let $m_1 \in \M_n$ such that 
\[ \ln(n)^4 \leq \frac{c^{\prime} d n }{c_{\mathrm{rich}} \ln(n)} \leq D_{m_1} \leq  \frac{c^{\prime} n }{\ln(n)} \leq n \enspace .\]
From \hypPrichMin, this is possible as soon as $n \geq L_{c_{\mathrm{rich}}, c^{\prime}, d}$. By \eqref{HP.eq.threshold} in Lemma~\ref{HP.le.threshold}, $A_n(m_0) \geq 1$ with probability at least $1 - L n^{-2}$. \\
We then have
\begin{align*}
\perte{\bayes_{m_1}} &\leq L_{\hypHPm,c^{\prime}}  \ln(n)^{\betamaj} n^{- \betamaj} \qquad \mbox{by \hypAp} \\
p_2(m_1) &\geq \paren{1 - \frac{L_{\hypHPm}}{\ln(n)} } \E\croch{p_2(m_1)}
\qquad \mbox{by \eqref{HP.eq.conc.p2.Qmp}} \\
\pen(m_1) &\leq K \E \croch{p_2(m_1)} \qquad \mbox{by \eqref{HP.eq.pen.mini}} \\
\absj{\delc(m_1)} &\leq L_A \sqrt{\frac{\ln(n)}{n}} \qquad \mbox{by \eqref{HP.eq.conc.delta.borne}}
\end{align*}
so that 
\begin{align*}
\crit^{\prime}(m_1) &\leq L_{\hypHPm,c^{\prime}}  \ln(n)^{\betamaj} n^{- \betamaj} + \paren{K - 1 + \frac{L_{\hypHPm}}{\ln(n)} } \E\croch{p_2(m_1)} + L_A \sqrt{\frac{\ln(n)}{n}} \\
&\leq \frac{(K - 1 + L_{\hypHPm} (\ln(n))^{-1}) \sigmin^2 c^{\prime}}{2 \ln(n)} 
\end{align*}
if $n \geq L_{\hypHPm,c^{\prime}}$.

We now choose $d$ such that the constant $d  L_{\hypHPm}$ appearing in the lower bound on $\crit^{\prime}(m)$ for ``small'' models is smaller than $(1 - K - L_{\hypHPm} (\ln(n))^{-1}) \sigmin^2 c^{\prime} /2$, that is $d \leq L_{\hypHPm,c^{\prime}}$. Then, we assume that $n \geq n_0 = L_{\hypHPm,c^{\prime},d}= L_{\hypHPm}$. Finally, we remove this condition as before by enlarging $\Kminpr$.
\end{enumerate}

\paragraph{Risk of $D_{\mh}$}
The proof of \eqref{HP.eq.mini.risque} is quite similar to the one of Lemma~\ref{HP.le.borne.Dmo}.
First, for every model $\mM_n$ such that $A_n(m) \geq 1$ and $D_m \geq \Kmindim n \ln(n)^{-1}$, we have 
\[ \perte{\ERM_m} \geq \punmin(m) \geq L_{\hypHPm} \Kmindim \ln(n)^{-2} \qquad  \mbox{ by \eqref{HP.eq.conc.p1.min.2}} \enspace .\]
Then, the model $m_0 \in \M_n$ defined previously satisfies $A_n(m) \geq 1$, and 
\[\perte{\ERM_{m_0}} \leq L_{\hypHPm} \paren{n^{- \betamaj / 2} + n^{-1/2}} \enspace .\]
If $n \geq L_{\hypHPm}$, the ratio between these two bounds is larger than $\ln(n)$, so that \eqref{HP.eq.mini.risque} holds. \BlackBox

\subsection{Concentration inequalities used in the main proofs} \label{HP.sec.conc}
In this section, we no longer assume that each model is the set of piecewise constant functions on some partition of $\X$. 
First, we control $\delc(m)$ with general models and bounded data.
\begin{proposition} \label{HP.pro.conc.delta.borne} 
Assume that $\norm{Y}_{\infty} \leq A < \infty$. Then for all $x \geq 0$, on an event of probability at least $1 - 2 e^{-x}$:
\begin{equation}\label{HP.eq.conc.delta.borne}
\forall \eta >0, \quad \absj{\delc(m)} \leq \eta \perte{\bayes_m} + \left( \frac{4}{\eta} + \frac{8}{3} \right) \frac{A^2 x}{n}  \enspace .\end{equation}
If moreover 
\begin{equation} \label{HP.def.Qmp}
\Qmp \egaldef \frac{n \E\croch{p_2(m)}}{D_m} > 0 \enspace , \end{equation}
on the same event,
\begin{align} \label{HP.eq.conc.delta.borne.2} \absj{\delc(m)} &\leq \frac{\perte{\bayes_m}}{\sqrt{D_m}} + \frac{20}{3} \frac{A^2 }{ \Qmp  }  \frac{\E[p_2(m)]}{\sqrt{D_m} } x 
\enspace .\end{align}
\end{proposition}
\begin{remark}[Regressogram case]
If $S_m$ is the set of piecewise constant functions on some partition $\paren{\Il}_{\lamm}$ of $\X$,
\[ \Qmp = \frac{1}{D_m} \sum_{\lamm} \sigl^2 \geq \paren{\sigmin}^2 >0\enspace .\]
\end{remark}

Then, we derive a concentration inequality for $p_2(m)$ in the regressogram case from a general result by \cite{Bou_Mas:2004}. 
\begin{proposition} \label{HP.pro.conc.p2}
Let $S_m$ be the model of piecewise constant functions associated with the partition $\paren{\Il}_{\lamm}$. Assume that $\norm{Y}_{\infty} \leq A$ and define $p_2(m) = P_n \paren{ \gamma(\bayes_m) - \gamma(\ERM_m) }$. 

Then, for every $x \geq 0$, there exists an event of probability at least $1-e^{1-x}$ on which for every $\theta \in (0;1)$,
\begin{equation} \label{HP.eq.conc.p2}
\absj{p_2(m) - \E\croch{p_2(m)}} \leq L \croch{\theta \perte{\bayes_m} + \frac{A^2 \sqrt{D_m} \sqrt{x}}{n} + \frac{A^2 x}{\theta n} }
\end{equation} 
for some absolute constant $L$.
If moreover $\sigma(X) \geq \sigmin >0$ a.s., we have on the same event:
\begin{equation} \label{HP.eq.conc.p2.Qmp}
\absj{p_2(m) - \E\croch{p_2(m)}} \leq \frac{L}{\sqrt{D_m}} \croch{\perte{\bayes_m} + \frac{A^2 \E\croch{p_2(m)} }{\sigmin^2} \paren{\sqrt{x} + x} } \enspace .
\end{equation}
\end{proposition}

Finally, we recall a concentration inequality for $p_1(m)$ proved by \cite[Proposition~9]{Arl:2008a}. Its proof is particular to the regressogram case. 
\begin{proposition}[Proposition~9, \cite{Arl:2008a}] \label{HP.pro.conc.p1}
Let $\gamma>0$ and $S_m$ be the model of piecewise constant functions associated with the partition $\paren{\Il}_{\lamm}$. Assume that $\norm{Y}_{\infty} \leq A < \infty$, $\sigma(X) \geq \sigmin >0$ a.s. and $\min_{\lamm} \set{n \pl} \geq B_n > 0 $.
Then, if $B_n \geq 1$, on an event of probability at least $1 - L n^{-\gamma}$, 
\begin{align} \label{HP.eq.conc.p1.min}
\punmin(m) &\geq \E \croch{\punmin(m)} - L_{A,\sigmin,\gamma}  \croch{  \frac{\ln(n)^{2}}{\sqrt{D_m}} +  e^{-L B_n} }  \E\croch{p_2(m)}
\\
\label{HP.eq.conc.p1.maj}
\punmin(m) &\leq \E \croch{\punmin(m)} + L_{A,\sigmin,\gamma}  \croch{  \frac{\ln(n)^{2}}{\sqrt{D_m}}  +  \sqrt{D_m} e^{-L B_n} }  \E\croch{p_2(m)} \enspace .
\end{align}
If we only have a lower bound $B_n > 0$, then, with probability at least $1-Ln^{-\gamma}$,
\begin{equation}
\label{HP.eq.conc.p1.min.2}
\punmin(m) \geq \paren{\frac{1}{2 + (\gamma+1) B_n^{-1} \ln(n) } -  \frac{L_{A,\sigmin,\gamma} \ln(n)^{2}}{\sqrt{D_m}} } \E\croch{p_2(m)} \enspace .
\end{equation}
\end{proposition}

\subsection{Additional results needed}
A crucial result in the proofs of Theorems~\ref{HP.the.opt.gal} and~\ref{HP.the.mini} is that $p_1(m)$ and $p_2(m)$ are close in expectation; the following proposition was proved by \citet[Lemma~7]{Arl:2008a}.
\begin{proposition}[Lemma~7, \cite{Arl:2008a}] \label{HP.pro.esp.p1p2}
Let $S_m$ be a model of piecewise constant functions adapted to some partition $\paren{\Il}_{\lamm}$. Assume that $\min_{\lamm} \set{n \pl} \geq B > 0$. Then, 
\begin{equation} 
\begin{split} \label{HP.eq.comp.Ep1.Ep2} 
\paren{1-e^{-B}}^2 \E \croch{p_2(m)} &\leq \E \croch{\punmin(m)} \\
&\leq \croch{ \mini{2}{\paren{1+ 5.1 \times B^{-1/4}}} + \paren{ \maxi{B}{1}} e^{- \paren{ \maxi{B}{1} } } } \E \croch{p_2(m)} \enspace . \end{split}
\end{equation}
\end{proposition}

Finally, we need the following technical lemma in the proof of the main theorems.
\begin{lemma}\label{HP.le.threshold}
Let $(\pl)_{\lamm}$ be non-negative real numbers of sum 1,
$(n\phl)_{\lamm}$ a multinomial vector of parameters
$(n;(\pl)_{\lamm})$. Then, for all $\gamma>0$, 
\begin{equation} 
\label{HP.eq.threshold} 
\min_{\lamm} \set{n \phl} \geq \frac{\min_{\lamm} \set{n \pl}}{2} - 2 (\gamma+1) \ln(n) \end{equation}
with probability at least $1 - 2 n^{-\gamma}$.
\end{lemma}
\begin{proof}
By Bernstein inequality \citep[Proposition~2.9]{Mas:2003:St-Flour}, for all $\lamm$, 
\[ \Prob \left( n \phl \geq (1-\theta) n\pl - \sqrt{2npx} - \frac{x}{3} \right) \geq 1 - e^{-x} \enspace .\]
Take $x = (\gamma+1) \ln(n)$ above, and remark that $\sqrt{2npx} \leq \frac{np}{2} + x$. The union bound gives the result since $\card(\Lambda_m) \leq n$.
\end{proof}

\subsection{Proof of Proposition~\ref{HP.pro.conc.delta.borne}} \label{HP.sec.proofs.delta}
Since $\norm{Y}_{\infty} \leq A$, we have $\norm{\bayes}_{\infty} \leq A$ and $\norm{\bayes_m}_{\infty} \leq A$. In fact, everything happens as if $S_m \cup \{ \bayes \}$ was bounded by $A$ in $L^{\infty}$.

We have 
\[ \delc(m) = \frac{1}{n} \sum_{i=1}^n \left( \gamma(\bayes_m,(X_i,Y_i)) - \gamma(\bayes,(X_i,Y_i)) - \E\left[ \gamma(\bayes_m,(X_i,Y_i)) - \gamma(\bayes,(X_i,Y_i)) \right] \right) \]
and assumptions of Bernstein inequality \citep[Proposition~2.9]{Mas:2003:St-Flour} are fulfilled with \[ c = \frac{8A^2}{3n} \quad \mbox{and} \quad v=\frac{8 A^2 \perte{\bayes_m}}{n} \] since
\[ \norm{ \gamma(\bayes_m,(X_i,Y_i)) - \gamma(\bayes,(X_i,Y_i)) - \E\left[ \gamma(\bayes_m,(X_i,Y_i)) - \gamma(\bayes,(X_i,Y_i)) \right] }_{\infty} \leq 8 A^2 \]
and \begin{align*} 
\var \left( \gamma(\bayes_m,(X_i,Y_i)) - \gamma(\bayes,(X_i,Y_i)) \right) &\leq \E \left[ \left( \gamma(\bayes_m,(X_i,Y_i)) - \gamma(\bayes,(X_i,Y_i)) \right)^2 \right] \\
&\leq 8 A^2 \perte{\bayes_m} 
\end{align*}
because $\norm{\bayes_m - \bayes}_{\infty} \leq 2A$ and
\begin{gather*}
\left( \gamma(t,\cdot) - \gamma(\bayes,\cdot) \right)^2 = \carre{ t(X) - \bayes(X)} \left( 2 (Y - \bayes(X)) - t(X) + \bayes(X) \right)^2  \\
\text{and} \quad \E \left[ (Y - \bayes(X))^2 \sachant X \right] \leq \frac{(2A)^2}{4} = A \enspace .
\end{gather*}
We obtain that, with probability at least $1 - 2 e^{-x}$, 
\begin{equation*}
\absj{\delc(m)} \leq \sqrt{2vx} + c = \sqrt{\frac{16 A^2 \perte{\bayes_m} x}{n}} + \frac{8 A^2 x}{3 n}
\end{equation*}
and \eqref{HP.eq.conc.delta.borne} follows since $2\sqrt{ab} \leq a \eta + b \eta^{-1}$ for all $\eta>0$. Taking $\eta = D_m^{-1/2} \leq 1$ and using $\Qmp$ defined by \eqref{HP.def.Qmp}, we deduce \eqref{HP.eq.conc.delta.borne.2}. \BlackBox

\subsection{Proof of Proposition~\ref{HP.pro.conc.p2}} \label{HP.sec.proof.p2}
We apply here a result by \citet[Theorem~2.2 in a preliminary version]{Bou_Mas:2004}, in which it is only assumed that $\gamma$ takes its values in $[0;1]$. This is satisfied when $\norm{Y}_{\infty} \leq A = 1/2$. When $A\neq 1/2$, we apply this result to $(2A)^{-1} Y$ and recover the general result by homogeneity. 

First, we recall this result in the bounded least-squares regression framework.
For every $t: \X \flens \R$ and $\epsilon >0$, we define \[ d^2(\bayes,t) = 2 \perte{t} \qquad \mbox{and} \qquad w(\epsilon) = \sqrt{2} \epsilon \enspace .\] 
Let $\phi_m$ belong to the class of nondecreasing and continuous functions $f: \R^+ \flens \R^+$ such that $x \flapp f(x)/x$ is nonincreasing on $(0;+\infty)$ and $f(1) \geq 1$. 
Assume that for every $u \in S_m$ and $\sigma>0$ such that $\phi_m(\sigma) \leq \sqrt{n} \sigma^2$,
\begin{equation} \label{HP.cond.phi}
\sqrt{n} \E \croch{ \sup_{t \in S_m, \, d(u,t) \leq \sigma} \absj{\perteempc{u} - \perteempc{t}} } \leq \phi_m(\sigma) \enspace . \end{equation}
Let $\esm$ be the unique positive solution of the equation
\begin{equation*}
\sqrt{n} \esm^2 = \phi_m(w(\esm)) \enspace . \end{equation*}
Then, there exists some absolute constant $L$ such that for every real number $q \geq 2$ one has
\begin{equation} \label{HP.eq.conc.p2.regression}
\norm{p_2(m) - \E[p_2(m)]}_q \leq \frac{L}{\sqrt{n}} \croch{ \sqrt{2 q} \paren{ \maxi{ \sqrt{\perte{\bayes_m}} } {\esm} }  + q \frac {2} {\sqrt{n}} } \enspace .\end{equation} 

\medskip

Using now that $S_m$ is the set of piecewise constant functions on some partition $\paren{\Il}_{\lamm}$ of $\X$, we can take 
\begin{equation} \label{HP.eq.phi_m.histos}
\phi_m(\sigma) = 3 \sqrt{2} \sqrt{D_m} \times \sigma \qquad \mbox{in \eqref{HP.cond.phi}.}\end{equation}
The proof of this statement is made below.
Then, $\esm = 6 \sqrt{D_m} n^{-1/2}$.

Combining \eqref{HP.eq.conc.p2.regression} with the classical link between moments and concentration 
\citep[see for instance][Lemma~8.9]{Arl:2007:phd}, the first result follows. The second result is obtained by taking $\theta = D_m^{-1/2}$, as in Proposition~\ref{HP.pro.conc.delta.borne}. \BlackBox

\paragraph{Proof of \eqref{HP.eq.phi_m.histos}}
Let $u \in S_m$ and $d(u,t) = \sqrt{2} \norm{u(X)-t(X)}_2$ for every $t: \X \flens \R$. 
Define $\psi: \R^+ \flens \R^+$ by
\[ \psi(\sigma) = \E \croch{ \sup_{d(u,t) \leq \sigma, \, t\in S_m}  \absj{ (P_n-P) (\gamma(u,\cdot) - \gamma(t,\cdot))} } \enspace .\]
We are looking for some nondecreasing and continuous function $\phi_m: \R^+ \flens \R^+$ such that $\phi_m(x)/x$ is nonincreasing, $\phi_m (1) \geq 1$ and for every $u \in S_m$, 
\[ \forall \sigma>0 \quad \mbox{such that} \quad \phi_m(\sigma)\leq \sqrt{n} \sigma^2 \enspace , \qquad \phi_m(\sigma) \geq \sqrt{n} \psi(\sigma) \enspace . \] We first look at a general upperbound on $\psi$.

\medskip

Assume that $u = \bayes_m$. If this is not the case, the triangular inequality shows that $\psi_{\text{general }u} \leq 2 \psi_{u = \bayes_m}$. Let us write
\[ t = \sum_{\lamm} \tl \1_{\Il} \qquad u = \bayes_m = \sum_{\lamm} \betl \1_{\Il} \enspace .\]

\paragraph{Computation of $P(\gamma(t,\cdot)- \gamma(\bayes_m,\cdot))$}
for some general $t \in S_m$:
\begin{align*}
P(\gamma(t,\cdot)- \gamma(\bayes_m,\cdot)) &= \E \left[ (t(X) - Y)^2 - (\bayes_m(X) - Y)^2 \right] \\
&= \E \left[ (t(X) - \bayes_m(X))^2 \right] + 2 \E \left[ (t(X) - \bayes_m(X)) (\bayes_m(X) - \bayes(X)) \right] \\
&= \E \left[ (t(X) - \bayes_m(X))^2 \right] \\
&= \sum_{\lamm} \pl (\tl - \betl)^2 \end{align*}
since for every $\lamm$, $\E \croch{ \bayes(X) \sachant X \in \Il } = \betl$.

\paragraph{Computation of $P_n(\gamma(t,\cdot)- \gamma(\bayes_m,\cdot))$}
for some general $t \in S_m$: with $\eta_i = Y_i - \bayes_m(X_i)$, we have
\begin{align*}
P_n(\gamma(t,\cdot)- \gamma(\bayes_m,\cdot)) &= \frac{1}{n} \sum_{i=1}^n \left[ (t(X_i) - Y_i)^2 - (u(X_i) - Y_i)^2 \right] \\
&= \frac{1}{n} \sum_{i=1}^n (t(X_i) - u(X_i))^2 - \frac{2}{n} \sum_{i=1}^n \left[ (t(X_i) - u(X_i)) \eta_i \right] \\
&= \frac{1}{n} \sum_{i=1}^n \sum_{\lamm} (\tl - \ul)^2 \1_{X_i \in \Il} - \frac{2}{n} \sum_{i=1}^n \sum_{\lamm} (\tl - \ul) \1_{X_i \in \Il} \eta_i \enspace .
\end{align*}

\paragraph{Back to $(P_n-P)$}
We sum the two inequalities above and use the triangular inequality:
\begin{align*}
\absj{ (P_n - P)(\gamma(t,\cdot) - \gamma(u,\cdot)) } &\leq 
\absj{\frac{1}{n} \sum_{i=1}^n \sum_{\lamm} (\tl - \ul)^2 (\1_{X_i \in \Il} - \pl)} \\ &\quad + \absj{ \frac{2}{n} \sum_{i=1}^n \sum_{\lamm} (\tl - \ul) \1_{X_i \in \Il} \eta_i }  \\
&\leq \frac{2A}{n} \sum_{\lamm} \left[ \left( \sqrt{\pl}\absj{\tl - \ul} \right) \frac{\absj{\sum_{i=1}^n   (\1_{X_i \in \Il} - \pl)} }{\sqrt{\pl}} \right] \\ &\quad +  \frac{2}{n}  \sum_{\lamm} \left[ \left(\sqrt{\pl}\absj{\tl - \ul} \right) \frac{\absj{ \sum_{i=1}^n \1_{X_i \in \Il} \eta_i }}{\sqrt{\pl}} \right]
\end{align*}
since $\absj{\tl - \ul} \leq 2A$ for every $t\in S_m$.

We now assume that $d(u,t) \leq \sigma$ for some $\sigma>0$, that is
\[ d(u,t)^2 = 2 \sum_{\lamm} \pl (\tl-\ul)^2 \leq \sigma^2 \enspace .\]
From Cauchy-Schwarz inequality, we obtain for every $t \in S_m$ such that $d(u,t) \leq \sigma$
\begin{align*}
\absj{ (P_n - P)(\gamma(t,\cdot) - \gamma(u,\cdot)) } &\leq 
\frac{2A\sigma} {\sqrt{2} n} \sqrt{ \sum_{\lamm} \frac{ \left(\sum_{i=1}^n   (\1_{X_i \in \Il} - \pl) \right)^2 }{\pl}} \\ 
&\quad + \frac{ \sqrt{2} \sigma}{n} \sqrt{ \sum_{\lamm} \frac{ \left(  \sum_{i=1}^n \1_{X_i \in \Il} \eta_i \right)^2} {\pl}}
\end{align*}

\paragraph{Back to $\psi$}
The upper bound above does not depend on $t$, so that the left-hand side of the inequality can be replaced by a supremum over $\set{t \in S_m \telque d(u,t) \leq \sigma}$. Taking expectations and using Jensen's inequality ($\sqrt{\cdot}$ being concave), we obtain an upper bound on $\psi$:
\begin{equation} \label{HP.eq.psi.intermediaire}
\begin{split}
\psi(\sigma) \leq \frac{2A\sigma} {\sqrt{2} n} \sqrt{ \sum_{\lamm} \E \left[ \frac{ \left(\sum_{i=1}^n   (\1_{X_i \in \Il} - \pl) \right)^2 }{\pl} \right]} 
 + \frac{\sqrt{2}\sigma}{n} \sqrt{ \sum_{\lamm} \E\left[ \frac{ \left(  \sum_{i=1}^n \1_{X_i \in \Il} \eta_i \right)^2} {\pl} \right]} 
 \end{split}
\end{equation}

For every $\lamm$, we have
\begin{equation} \label{HP.eq.psi.premier-terme} 
 \E \left(\sum_{i=1}^n   (\1_{X_i \in \Il} - \pl) \right)^2  = \sum_{i=1}^n \E \paren{ \1_{X_i \in \Il} - \pl }^2 = n \pl \paren{1-\pl} 
 \end{equation}
  which simplifies the first term. For the second term, notice that 
\begin{gather*} \forall i \neq j, \qquad \E \croch{ \1_{X_i \in \Il} \1_{X_j \in \Il} \eta_i \eta_j } = \E \croch{ \1_{X_i \in \Il} \eta_i } \E \croch{\1_{X_j \in \Il} \eta_j } \\
\mbox{and} \quad \forall i, \qquad \E \croch{ \1_{X_i \in \Il} \eta_i } = \E \croch{  \1_{X_i \in \Il} \E \croch{ \eta_i \sachant \1_{X_i \in \Il} } } = 0 \end{gather*} since $\eta_i$ is centered conditionally to $\1_{X_i \in \Il}$.
Then, 
\begin{equation} \label{HP.eq.psi.second-terme} 
\E \paren{\sum_{i=1}^n \1_{X_i \in \Il} \eta_i }^2 = \sum_{i=1}^n \E \croch{ \1_{X_i \in \Il} \eta_i^2 } \leq n \pl \norm{\eta}_{\infty}^2 \leq n \pl (2A)^2 \enspace . \end{equation}

Combining \eqref{HP.eq.psi.intermediaire} with \eqref{HP.eq.psi.premier-terme} and \eqref{HP.eq.psi.second-terme}, we deduce that
\begin{align*}
\psi(\sigma) \leq \frac{2A\sigma} {\sqrt{2} \sqrt{n}} \sqrt{D_m-1} + \frac{2\sqrt{2} A \sigma }{\sqrt{n}} \sqrt{ D_m} 
\leq 3 A \sqrt{2} \frac{\sqrt{D_m} } { \sqrt{n}} \times \sigma \enspace .\end{align*}
As already noticed, we have to multiply this bound by 2 so that it is valid for every $u \in S_m$ and not only $u = \bayes_m$.

The resulting upper bound (multiplied by $\sqrt{n}$) has all the desired properties for $\phi_m$ since $6 A \sqrt{2} \sqrt{D_m} = 3 \sqrt{2 D_m} \geq 1$. The result follows.
\BlackBox

\vskip 0.2in
\bibliography{arlot08a}

\end{document}